\documentclass[10pt,letterpaper]{amsart}
\usepackage{./macros}
 \usepackage{amsaddr}
 \usepackage{url}

\calclayout

\begin{document}


\title[Parallel, Adaptive Geometric Multigrid]{A Flexible, Parallel, Adaptive Geometric Multigrid Method for FEM}

\author{Thomas C. Clevenger}
\address{{clevenger.conrad@gmail.com}, Clemson University}
%
%
\author{Timo Heister}
\address{{heister@clemson.edu}, Clemson University}
%
\author{Guido Kanschat}
\address{kanschat@uni-heidelberg.de, Heidelberg University}
%
\author{Martin Kronbichler}
\address{kronbichler@lnm.mw.tum.de, Technical University of Munich}

\begin{abstract}
  Note: This is a preprint of the published article in ACM TOMS, see \url{https://doi.org/10.1145/3425193}.
  
  We present the design and implementation details of a geometric
  multigrid method on adaptively refined meshes for massively parallel
  computations. The method uses local smoothing on the refined part of
  the mesh. Partitioning is achieved by using a space filling curve
  for the leaf mesh and distributing ancestors in the hierarchy based
  on the leaves. We present a model of the efficiency of mesh
  hierarchy distribution and compare its predictions to runtime
  measurements.
  The algorithm is implemented as
  part of the \dealii{} finite element library and as such available to
  the public.
\end{abstract}

%
%
%





\maketitle

\section{Introduction}

Geometric multigrid methods are known to be solvers for elliptic
partial differential equations with optimal complexity in the number
of total variables~\cite{Hackbusch85,Bramble93}, but optimal performance
in a massively parallel environment depends on more than complexity
alone. Sufficiently many concurrent operations must allow utilization
of a sufficiently large part of the system, and it is not clear a
priori if multigrid methods with their hierarchy of coarse meshes and
synchronization due to grid transfer will be efficient on such
systems. In this article, we present algorithms
for such a method and demonstrate its feasibility in experiments.

Geometric multigrid methods for adaptive meshes and their
implementation on parallel computers have been studied for almost four
decades, for instance by~\cite{ZaveRheinboldt79,BankHolst03,Bastian96}
and others. A breakthrough was obtained in
the late 1990s by the use of space filling curves
(see~\cite{Zumbusch03} and literature cited therein), which allow the
partitioning of a hierarchical mesh in almost no time. Thus, load
balancing was reduced from an np-hard problem to a negligible task.
Such methods were implemented for instance in the software libraries
p4est~\cite{BursteddeWilcoxGhattas11},
\dealii{}~\cite{BangerthBursteddeHeisterKronbichler11}, DUNE~\cite{bastian08dune},
and Peano~\cite{WeinzierlMehl11,weinzierl2019peano}.

Several different kinds of adaptive multigrid methods can be
distinguished from the types of meshes and level spaces. Meshes can
either be conforming or nonconforming. Conforming meshes are generated by bisection, or by refinement into
$2^d$ children in $d$ dimensions dividing all edges and subsequent
closure (red-green refinement). These methods
are most commonly used with simplicial meshes, and typically require
mixed topologies otherwise.
The alternative are nonconforming methods,
most prominently the one-irregular meshes introduced by Bank, Sherman, and Weiser~\cite{BankShermanWeiser}. Here, the
difference in refinement between two cells sharing a common edge may
not exceed one level. This constraint is not
necessary and there have been codes which allow arbitrarily different
refinement levels of neighbors (see~\cite{aulisa2019construction} and references
therein). It nevertheless simplifies the code
considerably in particular in view of modern architectures. This
method has been implemented for simplicial meshes as well as meshes
based on (deformed) hypercubes. Since the meshes are nonconforming,
additional care has to be taken to ensure conformity of associated
finite element spaces. This is achieved by ``elimination of hanging
nodes'' resulting in algebraic constraints on the possible finite
element functions on the finer cell, see for instance \cite{ZaveRheinboldt79}.

After a locally refined mesh has been constructed, typically in an
adaptive algorithm, and its finite element space has been properly
defined with or without ``hanging nodes'', the resulting mesh has
cells on different levels. Thus, using a multigrid algorithm employing
smoothing operations on all cells on ``level $\ell$ or less'' is not
of optimal complexity on arbitrary meshes. Two remedies have been
proposed: local smoothing~\cite{Brandt77,Kanschat04,JanssenKanschat11,Chombo15}
and global coarsening~\cite{McCormickThomas86,BeckerBraack00,SundarBirosBursteddeRudiGhattasStadler12}.
We apply the former approach
in this work.
While optimal parallel work balance is more of a challenge with local smoothing,
there are a few potential advantages over global coarsening that justify
the investigation in this paper.
First, the computational complexity for the former is slightly lower and
optimal on all meshes, while there are (extreme) examples for
suboptimal complexity of global coarsening. Second, the smoothing
operation is always run on meshes without hanging nodes; while this is
not an issue for point smoothers like the Jacobi method, it
facilitates block smoothers, in particular patch smoothers as
in~\cite{ArnoldFalkWinther97hdiv,KanschatMao15}.
Finally, implementation
on vectorizing and multicore architectures is fairly straight-forward and does not
require special care at hanging nodes.

The hardware properties of state-of-the-art supercomputers have evolved
most rapidly in terms of the node-level
performance in the last decade, whereas network topologies across the nodes and node
numbers have been relatively steady with a thousand to ten thousand nodes on the top machines.
For these reasons, algorithmic components and data structures that
have low communication requirements are essential to balance inter-node latencies with
increasing intra-node performance, which can rely on hybrid parallelism, matrix-free
algorithms to relax the memory bandwidth requirements \cite{bauer17matrixfree,bauer18matrixfree,kronbichler2016comparison,kronbichler2019fast}, as
well as wide vectorization or offloading to GPUs, see e.g.~\cite{kronbichler2019gpu} and references therein.
These node-level optimizations provide fast
matrix-vector products for use in smoothers and level transfer,
including nearest-neighbor communication in the network.
The main focus of the present work
is on the algorithmic framework of local smoothing multigrid targeting
the inter-node case of large-scale parallel computations with MPI on meshes with adaptive
refinement. The components are flexible and allow for an
arbitrary element degree, various conforming or non-conforming elements,
as well as systems of equations, extending previous work on massively parallel multigrid \cite{MayBrown14,HPGMGv1,SundarStadlerBiros15,Rudi15,GholamiMalhotraSundarBiros16}.
Our contribution is integrated into the {\dealii{}} finite element library
and available as open-source software \cite{dealII91,dealii2019design}.

Our approach shown here can be summarized as a geometric multigrid method
on adaptively refined meshes. Each V-cycle is built from smoother, transfer, and coarse solver
operators that are equivalent to the serial method, while the work and data structures
are distributed in parallel. In practice, it is common to use Krylov methods,
like conjugate gradient, and perform a multigrid cycle as a preconditioner instead
of using multigrid as the solver.
This is the approach we use in the numerical examples.

The remainder of this work is structured as follows. In \S\ref{sec:gmg},
we present the geometric multigrid algorithm based on local smoothing. The
components for parallel execution in terms of the mesh infrastructure,
supported by an efficiency analysis of one particular partitioning strategy, are given in \S\ref{sec:parallel}.
Performance results are shown in \S~\ref{sec:performance} and the work is
concluded in \S~\ref{sec:conclusions}.

\section{Geometric multigrid with local smoothing}
\label{sec:gmg}
\subsection{Bilinear forms and finite element discretization}
The basis for our method is a partial differential equation in weak
form, abstractly written as: find $u\in V$ such that
\begin{gather}
  \label{eq:weak}
  a(u,v) = f(v) \qquad \forall v\in V.
\end{gather}
Here, $V$ is a suitable solution space and $f \in V^*$. For example, the Poisson
equation with homogeneous Dirichlet boundary condition on the domain
$\Omega \subset \mathbb R^d$ and right hand side $f\in L^2(\Omega)$
translates to $V=H^1_0(\Omega)$ and the weak equation
\begin{gather}
  a(u,v) \equiv \int_\Omega \nabla u\cdot \nabla v \,d x
  = \int_\Omega fv\,d x \equiv f(v).
\end{gather}

Our second example are the Lamé--Navier equations of linear elasticity
in space dimension $d$, where $V=H^1_0(\Omega;\mathbb R^d)$. With the
strain operator $\epsilon(u) = \frac12(\nabla u + \nabla u^T)$, we
obtain
\begin{gather}\label{eq:elasticity}
  a(u,v) \equiv \int_\Omega \Bigl[2\mu\epsilon(u): \epsilon(v)
  + \lambda \nabla\!\cdot u\;\nabla\!\cdot v\Bigr]\,d x
  = \int_\Omega f\cdot v\,d x \equiv f(v).
\end{gather}

These weak forms are discretized by the finite element
method. To this end, we introduce a mesh $\mesh_L$ covering the
domain $\Omega$. This article describes functionality of the library
{\dealii{}}, see~\cite{dealII91}, where the mesh cells $T$ are
quadrilaterals and hexahedra in two and three space dimensions,
respectively. We use mapped elements and the mappings from the
reference cell to the actual grid cell is not restricted to $d$-linear
functions or polynomials, but can be any function. On each mesh cell,
we define a local shape function space, typically by mapping
polynomials defined by a set of interpolation points from the
reference cell $[0,1]^d$. Using degrees of freedom, we establish
continuity between cells and define a basis of the finite element
space $V_L \subset V$ on the mesh $\mesh_L$. In the conforming case, the finite
element discretization of equation~\eqref{eq:weak} becomes: find
$u_L\in V_L$ such that
\begin{gather}
  \label{eq:conforming}
  a(u_L,v_L) = f(v_L) \qquad \forall v_L\in V_L.
\end{gather}
We will not distinguish between finite element functions $u_L\in V_L$
and their coefficient vectors $u_L\in \mathbb R^{n}$ with $n=\text{dim} V_L$,
since the meaning
will be clear from context. The basis used for this identification
consists of standard nodal finite element functions with local
support.

Discontinuous Galerkin (DG) finite element methods are an alternative
to conforming methods. Starting with the same mesh, we introduce
finite element spaces $V_L$ which are no longer conforming to the space
$V$, i.e. $V_L\not\subset V$, in particular spaces with no continuity
requirements. Therefore, the straight-forward discretization
using~\eqref{eq:conforming} is inconsistent and typically not
converging to the continuous solution. This is remedied by introducing
so-called \emph{flux terms} on the interfaces, which guarantee
consistency and stability of the method. Accordingly, the bilinear
form on $\mesh_L$ depends on the mesh itself and we write: find
$u_L\in V_L$ such that
\begin{gather}
  \label{eq:dg}
  a_L(u_L,v_L) = f_L(v_L) \qquad \forall v_L\in V_L.
\end{gather}
As an example, we mention the interior penalty method~\cite{Arnold82}
for the Laplacian with its multilevel analysis
in~\cite{GopalakrishnanKanschat03} and the bilinear form
\begin{multline}
  \label{eq:ip}
  a_L(u,v) \equiv \sum_{\cell\in\mesh_L}\int_\cell \nabla u\cdot\nabla v\,d x
  + \sum_{\face\in \facemesh_L^b} \int_F \Bigl[
  \sigma_h uv - \partial_n uv - u\partial_n v\Bigr]\,d s
  \\
  + \sum_{\face\in \facemesh_L^i} \int_F \Bigl[
  \sigma_h \jump{u}\jump{v}
  - 2\mean{\nabla u}\mean{v\mathbf n}
  - 2\mean{u\mathbf n}\mean{\nabla v}\Bigr]\,d s.
\end{multline}
The right-hand side with Dirichlet boundary data $u^D$ is
\begin{gather}
  f_L(v) \equiv \sum_{\cell\in\mesh_L} \int_\cell fv\,d x
  + \sum_{\face\in \facemesh_L^b} \int_F \Bigl[
  \sigma_h u^Dv - u^D\partial_n v\Bigr]\,d s.
 \end{gather}
Here, $\facemesh_L^i$ are the $(d-1)$-dimensional interfaces between
mesh cells of $\mesh_L$ and $\facemesh_L^b$ are the facets of cells on
the boundary of $\Omega$. Every face $\face\in \facemesh_L^i$ has two
adjacent cells, say $\cell^+$ and $\cell^-$. We call the restriction
of the finite element functions $u$ and $v$ to these cells $u^+$,
$u^-$, $v^+$, and $v^-$, respectively. With these definitions, we have
the jump and mean value operators
\begin{gather}
    \jump{u} = u^+-u^-,\qquad
    \mean{u} =\frac{u^++u^-}{2}.
\end{gather}
Finally, $\sigma_h$ is an appropriately chosen penalty parameter inversely
proportional to the diameter of the cells attached to the face $F$.

\subsection{Geometric multigrid}
\label{sec:gmg-serial}

The geometric multigrid method employs a hierarchy of meshes, which
we generate as follows. Starting from the coarse mesh $\mesh_0$,
we generate the mesh
 $\mesh_{\ell+1}$ from $\mesh_\ell$ by selecting a subset or all of its
cells and refining these isotropically by bisecting each edge,
generating $2^d$ children each.
This results in a sequence
\begin{gather}
  \label{eq:mesh-hierarchy}
  \mesh_0 \sqsubset \mesh_1 \sqsubset \dots \sqsubset \mesh_L,
\end{gather}
where the symbol ``$\sqsubset$'' denotes nested meshes, that is, every
cell of a mesh on the left of this symbol is the union of one or more
cells of the mesh on the right.

 As usual, we define finite element
spaces $V_\ell$ on these meshes by defining local shape function
spaces on each cell $\cell\in \mesh_\ell$ and concatenating these
spaces, identifying shape functions on adjacent cells which are
associated to joint degrees of freedom. For most finite elements, and
these are the ones we consider here, the shape functions on a cell
$\cell$ can be represented as linear combinations of the shape
functions on its children in the mesh hierarchy. Therefore, the mesh
hierarchy above induces a sequence of finite element spaces
\begin{gather}
  \label{eq:space-hierarchy}
  V_0 \subset V_1 \subset \dots \subset V_L.
\end{gather}

We discretize the weak formulation~\eqref{eq:weak} on each mesh by a
bilinear form $a_\ell(.,.)$ and the problem: find $u_\ell \in V_\ell$,
such that
\begin{gather}
  a_\ell(u_\ell, v_\ell) = f_\ell(v_\ell)\quad\forall v_\ell\in V_\ell.
\end{gather}
For conforming finite element methods, the bilinear forms and the
right hand side are simply the restrictions of $a(.,.)$ and $f(.)$ to
the space $V_\ell$. For DG and other stabilized schemes, they contain
additional terms for consistency and stability. Associated with the
bilinear form $a_\ell(.,.)$ is a linear operator
$A_\ell:V_\ell \to V_\ell$ defined by
\begin{gather}
  \left(A_\ell u_\ell, v_\ell\right)_{V_\ell}
  = a_\ell(u_\ell, v_\ell)
  \quad\forall u_\ell, v_\ell \in V_\ell.
\end{gather}
Here, the inner product on $V_\ell$ is the one used in the conjugate
gradient method, typically the Euclidean norm of the coefficient
vector of a function $u_\ell\in V_\ell$ with respect to the nodal
basis of $V_\ell$, see for instance the discussion
of mesh dependent norms in~\cite{BraessHackbusch83,BrennerScott08} and
their relation to the inner product of $L^2(\Omega)$. Based on the
embeddings in~\eqref{eq:space-hierarchy}, we define the grid transfer
operators
\begin{xalignat}2
\prol_\ell&:V_\ell\to V_{\ell+1}
& v&\mapsto v, \\
\rest_\ell&:V_{\ell+1} \to V_\ell
& \left(\rest_\ell u, v_{\ell}\right)_{V_\ell}
&= \left( u, v_{\ell}\right)_{V_{\ell+1}}
\;\forall v_\ell\in V_\ell.
\end{xalignat}

On each mesh level $\ell$, we employ a smoother
$\smoother_\ell(u_\ell, g_\ell)$, which employs the right hand side
$g_\ell$ and the current state $u_\ell$ to compute a result.
Examples for such smoothers are relaxation methods of the form
in Algorithm~\ref{alg:smoother}.
\lstset{caption={The abstract smoother algorithm.},label={alg:smoother}}
\begin{algorithm}
  function $\smoother_\ell$($u_\ell, g_\ell$)
   for $k:=1,\dots,m_\ell$
     $u_\ell \;\gets\; u_\ell + B_\ell^{-1}\bigl( g_\ell  - A_\ell u_\ell\bigr)$
   return $u_\ell$
\end{algorithm}
Here, $m_\ell$ is the number of smoothing steps and
$B_\ell$ is the type of relaxation method, for instance the diagonal
for the Jacobi method or the lower triangle for
Gauss--Seidel. Similarly, additive and multiplicative Schwarz methods
fit into this concept, but it also extends to nonlinear methods like
conjugate gradients or GMRES.

\lstset{caption={The abstract multigrid V-cycle algorithm.},label={alg:vcycle}}
\begin{algorithm}
  function Vcycle ($\ell, g_\ell$)
  if $\ell > 0$
    $u_1$ $\gets$ $\smoother_\ell$($0, g_\ell$)
    $u_2$ $\gets$ $u_1 + \prol_{\ell-1}$Vcycle$\left(\ell-1, \rest_{\ell-1}\left(g_\ell - A_\ell u_1\right)\right)$
    $u_3$ $\gets$ $\smoother_\ell$($u_2, g_\ell$)
    return $u_3$
  else
    return $A_0^{-1} g_\ell$
\end{algorithm}

We are now ready to state the multigrid V-cycle algorithm in Algorithm~\ref{alg:vcycle} in abstract
form, as it has been done in numerous publications.
In addition to the level transfers and smoothers discussed before, the
recursion of the algorithm requires closure at level 0, denoted as the
inverse of $A_0$. This is called coarse grid solver, and in an
implementation can be a direct solver since the system is small, or a
basic iterative method like conjugate gradients or GMRES since the
system is well conditioned.

\subsection{Level meshes and local smoothing}
For uniformly refined meshes, the definition of the level meshes $\mesh_\ell$ is trivial.
For adaptively refined meshes, we consider the hierarchy as a tree or a forest
(for more than one coarse cell).
Each node of this forest is
a cell in the mesh hierarchy. The cells of the coarse mesh $\mesh_0$ are the
roots of its trees.
The level of a cell is defined as the distance from its root.
The mesh on which we discretize the
differential equation consists of the leaves of this tree (or forest)
and is denoted as the \emph{leaf mesh} $\mesh_L$. Since it is
obtained by local refinement, it typically consists of cells on
different levels up to level $L$.

For intermediate levels $0<\ell<L$,
two different
schemes have been devised.
First, for global coarsening, the meshes are constructed from the leaf mesh. The mesh $\mesh_{\ell-1}$
is obtained by replacing all cells in $\mesh_{\ell}$ by their immediate parents
in the forest. Once a cell of the coarse mesh is reached, it remains in further level meshes.
For more information we refer
to~\cite{BeckerBraack00,SundarBirosBursteddeRudiGhattasStadler12}.

Instead, in our definition of the
level mesh
$\mesh_\ell$ consists of all cells of level $\ell$ and of all leaves
with level less than $\ell$.
With such a definition of $\mesh_\ell$, a fairly coarse leaf cell can be
part of many different level meshes.

\begin{figure}[tp]
  \centering
  \begin{tabular}{ccc}
    \includegraphics[width=.3\textwidth]{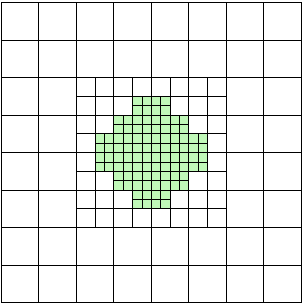}
    &
    \includegraphics[width=.3\textwidth]{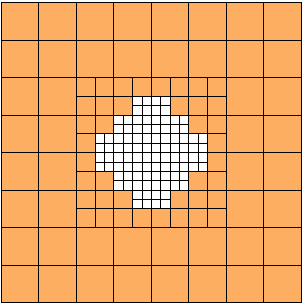}
    &
    \includegraphics[width=.3\textwidth]{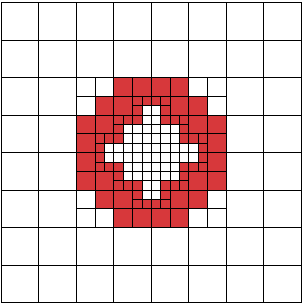}
    \\
    $\operatorname{supp} V_\ell^S = \mesh_\ell^S$
    &
    $\operatorname{supp} V_\ell^L = \mesh_\ell^L$
    &
    $\operatorname{supp} V_\ell^I$
  \end{tabular}
  \caption{Splitting of the mesh $\mesh_\ell$ and the space $V_\ell$
    into subspaces for local smoothing. Superscript $S$ refers to the
    cells and functions strictly on level $\ell$ (left), used for
    smoothing. Superscript $L$ is the support of functions actually
    defined on lower levels (center) and $I$ is the support of functions for node
    functionals on the interface which have support in both
    subdomains}
  \label{fig:level-splitting}
\end{figure}

In order to obtain an algorithm
with optimal complexity, smoothing for the degrees of freedom of a
given cell should only happen on a single level. This is where local
smoothing enters: while we are running a multigrid method for
the whole finite element space $V_\ell$, we restrict smoothing only to
the mesh cells which are actually on level $\ell$. This splitting is
explained in Figure~\ref{fig:level-splitting}.
The mesh $\mesh_\ell$ is split into the submesh $\mesh_\ell^S$ of
cells strictly on level $\ell$ and $\mesh_\ell^L$ of cells on lower
levels than $\ell$. For DG methods, this immediately results in a
splitting $V_\ell = V_\ell^S\oplus V_\ell^L$, where the support of
each subspace is its corresponding submesh. The splitting for
continuous methods is more complicated since there are finite element
basis functions with support straddling the interface and thus in both
$\mesh_\ell^S$ and $\mesh_\ell^L$. The span of these basis functions
is called $V_\ell^I$, the space of interface functions.

We now give a short review of the structure of the operators in the
multigrid method outlined
in~\cite{Kanschat04,JanssenKanschat11}. Here, the goal is to implement
the algebraic equivalent of the original multigrid method for the
space hierarchy $\{V_\ell\}$ with operators obeying the subspace
splitting.  We start with the observation that conforming methods
require the function on the refined side of a refinement edge to
coincide with the function on the coarse side. This translates into
elimination of degrees of freedom on the refined side such that
$V_\ell^I$ uses degrees of freedom representable on the coarse mesh
$\mesh_{\ell-1}$ only. The fact that shape functions in $V_\ell^I$
have different support from their counterparts in $v_{\ell-1}$ is
taken care of by the grid transfer operators.

Thus, we can restrict smoothing on
level $\ell$ to $V_\ell^S$ and can ignore $V_\ell^I$. Furthermore, in
the case of DG methods, $V_\ell^I = \{0\}$, such that in both cases we
can write $V_\ell = V_\ell^S \oplus V_\ell^I \oplus V_\ell^L$. Our assumptions on local smoothing translate to
\begin{gather}
  \smoother_\ell\left(
    \begin{pmatrix}
      x^S\\x^I\\x^L
    \end{pmatrix},
    \begin{pmatrix}
      g^S\\g^I\\g^L
    \end{pmatrix}\right)
    =
    \begin{pmatrix}
      \smoother_\ell^S(x^S,g^S)\\x^I\\x^L
    \end{pmatrix},
\end{gather}
where $\smoother_\ell^S(x^S,g^S)$ is now the local smoother on
$V_\ell^S$ only. We observe that the embedding operator $R_{\ell-1}^T$
maps a function from $V_{\ell-1}$ to itself.
Therefore, $R_{\ell-1}$ is the identity on
$V_\ell^L$. Thus, $R_{\ell-1}$ has the structure
\begin{gather}
  R_{\ell-1}
  \begin{pmatrix}
    x^S\\x^I\\x^L
  \end{pmatrix}
  = R_{\ell-1}^Sx^S + R_{\ell-1}^I x^I + x^L.
\end{gather}
Note that the restriction operator involves more than two levels if its range is not in $V_{\ell-1}^S$, because the other components of $V_{\ell-1}$ are not stored on this level. Therefore, we simplify the code considerably if we ensure $R_{\ell-1}^SV_\ell^S + R_{\ell-1}^IV_\ell^I \subset
V_{\ell-1}^S$. This can be achieved by an additional mesh refinement rule: require that there is a complete layer of cells on level
$\ell-1$ surrounding $\mesh_\ell^S$.

Residuals on the other hand must be computed
correctly on the whole space $V_\ell$ according to
\begin{gather}
  \begin{pmatrix}
    r_\ell^S\\r_\ell^I\\r_\ell^L
  \end{pmatrix}
  =
  \begin{pmatrix}
    g^S\\g^I\\g^L
  \end{pmatrix}
  -
  \begin{pmatrix}
    A_\ell^{S} & A_\ell^{SI} & A_\ell^{SL} \\
    A_\ell^{IS} & A_\ell^{I} & A_\ell^{IL} \\
    A_\ell^{LS} & A_\ell^{LI} & A_\ell^{L}
  \end{pmatrix}
  \begin{pmatrix}
    x^S\\x^I\\x^L
  \end{pmatrix}.
\end{gather}
Note that the matrices $A_\ell^{SL}$ and $A_\ell^{LS}$ are the flux
matrices of a DG method on the refinement edge and thus vanish for
conforming methods. Furthermore, we see in the V-cycle algorithm that
this residual is immediately restricted to the coarse space
$V_{\ell-1}$. Since the restriction acts as identity on $V_{\ell}^L$,
we can avoid computing $r_\ell^L$ and defer it to the lower
level. Thus, the matrix $A_\ell^L$ is not needed in computations at
all. The matrix $A_\ell^S$ is used for smoothing on level $\ell$. The
off-diagonal matrices correspond to coupling between degrees of
freedom on the cells at the interface, and are needed in addition to
$A_\ell^S$ for a consistent multigrid method.

A major advantage of local smoothing is its fairly simple data
structure. The level meshes $\mesh_\ell^S$ do not have hanging nodes,
such that the results of cell-wise operations can be entered into global
vectors very efficiently without any elimination process.
This also simplifies implementation, especially when considering more
involved smoothers such as patch smoothers, that would need to operate
on patches of cells on different refinement levels.
Furthermore,
it is of optimal computational complexity on any locally refined mesh,
while global coarsening may be suboptimal on some meshes with extreme
local refinement, see~\cite{JanssenKanschat11}. Nevertheless, this
second aspect does not seem to have much impact on actual
computations.

\section{Parallelization of geometric multigrid}
\label{sec:parallel}

We will now discuss the construction of an efficient and scalable
parallel version of the adaptive multigrid method described in
Section~\ref{sec:gmg}. We emphasize data and communication structures
while keeping the algorithm mathematically equivalent to the weathered sequential version. Regarding parallelism, we have to consider three levels of
parallelization in modern computer architectures, namely message
passing between computer nodes and intra-node parallelization
separated in multicore/multitasking (multiple instruction, multiple
data) and vectorization (single instruction, multiple data).
As motivated in the introduction, this article focuses on message passing. The intra-node parallelization approach employed is shortly discussed in \S\ref{sec:matrix-free}.

A scalable approach requires distributed data structures and scalable algorithms
operating on them including equal partitioning of the work.
As demonstrated in the computations in the later sections of this paper, the parallel algorithms described here enable high resolution adaptive computations with billions of
unknowns on 100,000+ cores.
We concentrate on MPI as the parallelization framework and we refer to a single
MPI rank or process as ``processor''.

\subsection{Parallel algorithm}

Our
algorithm is synchronized between applications of residual,
smoothing, grid transfer operators, and coarse grid solvers.  Hence,
our focus lies in the parallel implementation of these operators.

The abstraction of parallel data structures and algorithms equivalent to the serial version is well-known.
Libraries such as PETSc \cite{petsc-user-ref,petsc-web-page} and Trilinos \cite{trilinos-overview} have provided
linear algebra data structures (vectors, sparse matrices) and algorithms (iterative
solvers) with this abstraction for a long time. Up to a point, this isolates the user (for example finite element library implementors) from having to interface directly with the underlying parallel computing framework. Abstraction of this kind are of course not perfect, because
operations like finite element assembly need to be partitioned between the processors.
Nevertheless, it enables the design of parallel algorithms on a higher level, like
it is done in \dealii{}, see \cite{BangerthBursteddeHeisterKronbichler11}.

The workload is typically distributed by partitioning the cells of the computation using
graph based partitioners or using space-filling curves (like METIS \cite{metis}, Zoltan \cite{ZoltanIsorropiaOverview2012}, or p4est \cite{BursteddeWilcoxGhattas11} -- the latter one being used in \dealii{}). This partitioning can be used to distribute cell-based work,
like matrix or residual assembly, and can be used to generate a partitioning of degrees of freedom that
is needed for the row-wise division of linear algebra objects (vectors, matrices). The latter step requires
a rule to decide on the ownership of degrees of freedom on the interface between processor boundaries of the
cells.
For the finite element framework, the main effort is to correctly assign and communicate ghost cells and ghost indices, while
the communication for matrix-vector products and finite element assembly of foreign entities only involves
neighboring processors and is
typically provided by the linear algebra libraries.

Here, we will follow the same approach for the partitioning of cells and degrees of freedom on
each level of the multigrid hierarchy: after partitioning of all cells strictly on level $\ell$ in some way, we use this to partition the degrees of freedom accordingly.
Like above, it is advantageous for large computations if only part of the mesh that is relevant for
the current processor are stored locally.
There are different options for partitioning cells on each multigrid
level. We
will discuss different strategies and the approach we take in
Section~\ref{sec:partition-first-child}, but stress that our implementation is flexible in this respect.

While knowledge about the whole mesh is not required, we need ghost neighbors on each level, which can be on different levels in adaptive computations.
In our implementation, we decided to always construct the ghost neighbors as the set of all cells that share
at least a vertex with the local cells and exchange all information about them, even though specific implementations might require less information (for example only neighbors across faces for DG).
Furthermore, information about parents/siblings is required for transfer operations. This allows us to compute and exchange
the necessary information about degrees of freedom for smoothing and grid transfer.

To summarize the execution of multigrid in parallel, the following parallel ingredients are necessary:

\begin{itemize}
\item
Prolongation and restriction are conceptually a multiplication
of distributed vectors with a rectangular transfer matrix
and as such equivalent to the serial transfer. Known algorithms for sparse matrix-vector products
scale well in parallel.

\item
The smoothers we consider are conceptually a sequence of local operations on individual degrees of freedom
or small subspaces. Additive smoothers (Jacobi, additive Schwarz, etc.) can be run in parallel
on all processors and are still equivalent to the serial method. Sequential
smoothers (Gauss--Seidel, multiplicative Schwarz, etc.) can of course not be used immediately.

Other parallelizable smoothers (for example block Gauss--Seidel variants) may not be
equivalent to their sequential version. Only in this case equivalence to the
sequential multigrid algorithm is lost, which can lead to an increase in
iteration numbers when increasing the number of processors, albeit they often
remain more effective than additive smoothers.

\item
There are several options for coarse solvers. First, if the problem is reduced
to a very small number of cells and processors, runtime is negligible and
(parallel) direct solvers can be applied. In other cases, when the coarse mesh
still has a large number of degrees of freedom, switching to algebraic multigrid
is an option used in \cite{SundarBirosBursteddeRudiGhattasStadler12,SundarStadlerBiros15}
for example and what we use in Section~\ref{sec:test-elasticity}.

\end{itemize}

\subsection{Matrix-free implementation}
\label{sec:matrix-free}
In the previous sections, we have derived the multigrid algorithm in an abstract
way based on linear algebra operators. While these are typically implemented
as sparse matrices, the concept directly translates to matrix-free operator
evaluation. These methods often provide considerably faster evaluation of
matrix-vector products than assembled matrices, in particular for higher order
finite elements, because the access to memory is significantly
reduced~\cite{kronbichler2012generic,kronbichler2019fast}, which is the limiting factor in
matrix-based implementations. In this work, we consider methods based on
sum factorization techniques on hexahedra which have a particularly high
node-level performance \cite{kronbichler2016comparison} and are also applicable to
GPUs~\cite{kronbichler2019gpu}. We note that a fast intra-node performance puts more emphasis on possible
communication bottlenecks.

\subsection{Partitioning strategy for mesh hierarchy}
\label{sec:partition-first-child}

When partitioning the cells on each level of the multigrid hierarchy, a balance between several conflicting
goals is necessary:
\begin{enumerate}
\item Minimize communication for transfer operations between multigrid levels.
\item Fair work balance on each level (same number of cells per processor).
\item Minimize interface between processor boundaries on each level (minimizes communication in smoother applications or residual computations).
\item Minimize required additional storage for the mesh hierarchy if local cells have little overlap between levels.
\end{enumerate}

Previous work \cite{SundarBirosBursteddeRudiGhattasStadler12} has concentrated
mainly on aspect (2) by using an independent partition on each level. The
multigrid method there is based on global coarsening instead of local
smoothing, so each level is an adaptively refined mesh that needs to be
partitioned. While satisfying (2), this ignores (1) and requires duplicate
storage, violating (4). Experiments in
\cite{SundarBirosBursteddeRudiGhattasStadler12} have shown excellent scaling
despite these deficiencies. Note that (1) is satisfied for meshes refined
mostly globally, but duplicate storage (4) is still required. It is worth
pointing out that the storage requirement reduces by a factor of 8 in each coarsening
step in 3d.

In this work, we explore algorithms that ignore (2) and partition the hierarchy
based on the partitioning of the leaf mesh to minimize communication cost
and storage requirements (goal (1) and (4)). This option also represents a
straight-forward implementation in terms of local tree data structures and can
be seen as a baseline to more sophisticated setups.
We will see that we satisfy goal (2) for mostly globally refined meshes and that we can quantify the partitioning efficiency (see Section~\ref{sec:partition-efficiency}).

Note that both approaches behave similarly for uniformly refined meshes, while goals (1) and (2) are conflicting for an adaptive scheme. Finally, note that (2) is desired when assuming that levels are passed through sequentially, as the multigrid algorithm suggests, but one could design
a parallel method that does not require synchronization on each level.

In the following, we will partition the multigrid cells by the {\bf ``first-child rule''} as follows:
First, distribute the leaf cells using a space filling curve (we use p4est~\cite{BursteddeWilcoxGhattas11} as described in \cite{BangerthBursteddeHeisterKronbichler11}).
Second, for each cell in the hierarchy, recursively assign the parent of a cell to the owner of the first child cell.
A similar strategy has been used also by the Octor package \cite{Tu2005Octor} for the parent-child relationships.

For an example with seven cells, see Figure~\ref{fig:mesh7-partition} that shows the mesh with the
space-filling curve on the left, the tree representing the refinement in the middle,
and the cells on each level on the right.
This approach has the following consequences:
\begin{enumerate}
 \item The local cells and their parents are already present on each processor and the ownership of
 parents is known without any communication. This means the partitioning of the multigrid
 hierarchy can be done without communication.
 \item No duplicated storage for the mesh is needed as all parent cells are already stored locally (goal (4)).
 \item Transfer operations are local and require only a small amount of communication at processor boundaries
 (goal (1)).
 \item Processors drop out automatically on coarser levels, which is desired.
 \item The workload on each level is not distributed equally.
\end{enumerate}
We will discuss the last point and its impact in the next subsection.

\subsection{Partitioning efficiency model}
\label{sec:partition-efficiency}

Our model for the complexity of the partitioned workload, in short
parallel complexity, is based on the assumption that parallelization
is completely achieved by MPI ranks and that within each rank the workload is proportional to the number of cells. Below, we develop a complexity model based
on this assumption, estimating the parallel complexity of our
algorithm in terms of mesh cells per level.

Let $N_{\ell}$ be the number of cells on level $\ell$ and $N_{\ell,p}$ of the subset owned by processor
$p$. Here, cell refers to all cells in the hierarchy, not only leaf cells belonging to the discretization mesh.
We assume that the workload for each cell is equal, such that
$N_{\ell,p}$ is proportional to the total amount of work a processor
has to invest on level $\ell$. Obviously, the optimal
parallel complexity is
\begin{gather*}
  W_{\text{opt}} =  \frac1{n_{p}} \left[ W_0 +  \sum_{\ell=1}^L \sum_p N_{\ell,p} \right]
  = \frac1{n_{p}} \left[W_0 +  \sum_\ell N_{\ell}\right].
\end{gather*}
Here, the terms in brackets specify the total work of the multigrid algorithm and
the total number of processors is given by $n_p$.
$W_0$ is the cost of the coarse grid solver, which may
be different than the cost of a smoother application.

This calculation is based on perfect equidistribution of work and
neglects communication overhead. In particular, it is not achievable
if grid transfers are synchronized, as in our implementation. In this
case, we can only distribute the work on each level such that we are
bound from below on each level by
\begin{gather*}
  W_{\ell,\text{opt}} = \left\lceil \frac1{n_{p}}  N_{\ell} \right\rceil,
\end{gather*}
where $\lceil n \rceil$ is the smallest integer greater or equal to $n$.
Therefore, the best achievable work time with syncing between levels is
\begin{gather*}
  W_{\text{sync}} = W_0 + \sum_{\ell=1}^L W_{\ell,\text{opt}}.
\end{gather*}
On the other hand, with imperfect distribution of work, the limiting effort on each level is
\begin{gather*}
  W_\ell = \max_{p} N_{\ell,p},
\end{gather*}
and the total parallel complexity $W$ and partitioning efficiency $\parteff$ due to imbalance against a hypothetical optimal partitioning are given by
\begin{gather}
  W = W_0 + \sum_{\ell=1}^L W_\ell,
  \qquad
  \label{eq:slowdown-model}
  \parteff = \frac{W_{\text{opt}}}{W}.
\end{gather}
To summarize, the efficiency $\parteff$ quantifies the overhead produced by
not having a perfect work balance on each level of the multigrid
hierarchy. The work is modeled as one unit per cell and, as such, can be used
to model the computational cost for grid transfer, smoother application, and
other operations.
The efficiency, by definition, accumulates the imbalances on each level.

\begin{figure}[tp]
\tikzsetnextfilename{tikz-figures/tree7-grid-active}
\includegraphics[width=0.3\textwidth]{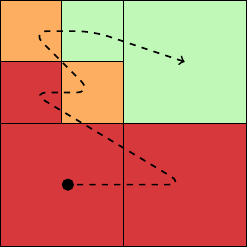}
\hfill
\tikzsetnextfilename{tikz-figures/tree7-notopt}
\includegraphics[height=4cm]{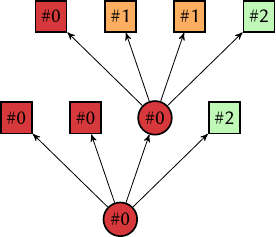}
\hfill
\tikzsetnextfilename{tikz-figures/tree7-grid-level}
\includegraphics[height=4cm]{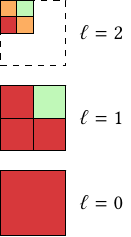}
 \caption{Example partitioning for an adaptive mesh from a single coarse cell with 7 leaf cells distributed between three processors indicated by color and label (red \#0, orange \#1, light green \#2). Ownership of the level cells
 is determined by the ``first-child rule'' (round shapes in the middle).
 Left: The partitioning of the leaf mesh with the space-filling $z$-curve.
 Middle: Tree representation of the mesh and its partitioning.
 Right: Partitioning of the three multigrid layers $\mesh_\ell^S$.
 }
  \label{fig:mesh7-partition}
\end{figure}

We give an example for these estimates for the mesh hierarchy displayed in Figure~\ref{fig:mesh7-partition}.
It consists of 7 leaf cells obtained by successive refinement of a single coarse
cell. The partitioning is done for three processors. The ownership of the leaf cells is determined by p4est using a
space-filling curve ($z$-curve, also known as Morton curve, dashed line on the left picture) or depth-first traversal (from left to right) in the tree
representation depicted in the middle. The ownership of cells in the multigrid hierarchy (round circles in the tree) is determined by copying the leaf ownership and then
applying the ``first-child rule'' recursively. For example, the parent of the four smallest cells on level 2 is red (\#0)
because the first (bottom-left) child also belongs to processor \#0 (red).
One result of this partitioning is that processors drop out on coarser
levels automatically. Here, processor \#1 (green) recuses itself on
level 1 and only processor \#0 (red) remains on the coarsest level
(here a single cell). The coarsest mesh is not necessarily completely
owned by processor \#0 if it consists of more than a single cell.

The optimal parallel complexity is simply the number of all cells in the hierarchy
divided by the number of processors, hence $W_{\text{opt}} = \frac{9}{3} = 3$. On
the other hand, assuming the coarse grid solver has the same
complexity as the work load per cell on higher levels,
\begin{gather*}
  W_{\text{sync}} = W_{0,\text{opt}} + W_{1,\text{opt}} + W_{2,\text{opt}}
  = \ceil{\tfrac{1}{3}} + \ceil{\tfrac{4}{3}} + \ceil{\tfrac{4}{3}} = 5.
\end{gather*}
Comparing to
\begin{gather*}
  W = \sum_\ell W_\ell = 1 + 2 + 3 = 6,
\end{gather*}
we obtain $\parteff= 1/2$ in this example.  In other words, our model predicts a slowdown of 100\%
and 20\% compared to $W_\text{opt}$ and $W_{\text{sync}}$, respectively.
The slowdown with respect to $W_{\text{sync}}$ is due to the
non-optimal partitioning on level 1, where processor \# 0 (red) works
on three cells while the other processors have to wait. An optimal
partitioning would only require operating on two cells sequentially on
that level. Compared to $W_{\text{opt}}$, we do not have enough
cells to keep three processors busy.

This example suggests that the efficiency of the algorithm depends
significantly on the base of comparison, $W_{\text{sync}}$ or
$W_{\text{opt}}$. In fact, a closer inspection of the definitions
reveals that they only differ by rounding up the load on each level to
the next multiple of $n_p$, a difference which drops below 1\% as soon
as we have 100 cells on each processor. Below, we only use
$W_{\text{opt}}$ when we assess the efficiency of our mesh hierarchy
distribution.

\subsection{Experimental study of the efficiency of the first-child rule}

\begin{figure}[tp]
  \centering
  \subfigure[``circle'' refinement. Leaf mesh with finest cells on level 9.]
  {
    \begin{tabular}{ccc}
        \includegraphics[width=.3\textwidth]{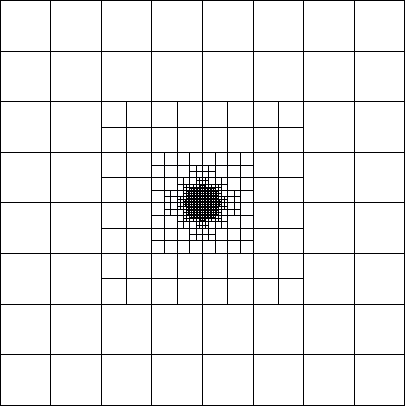}
    \end{tabular}
    \label{fig:circle-ref}
  }
  \hspace{2cm}
  \subfigure[``quadrant'' refinement. Leaf mesh with finest cells on level 8.]
  {
    \begin{tabular}{ccc}
        \includegraphics[width=.3\textwidth]{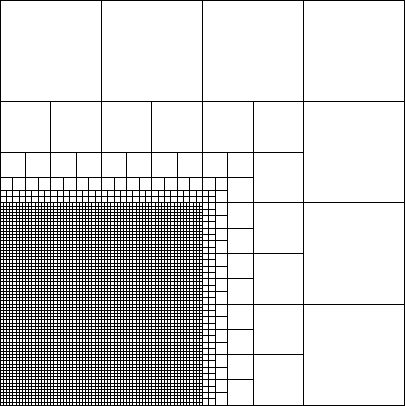}
    \end{tabular}
    \label{fig:quad-ref}
  }
  \subfigure[``sine curve'' refinement. Leaf mesh with finest cells on level 9.]
  {
    \begin{tabular}{ccc}
        \includegraphics[width=.3\textwidth]{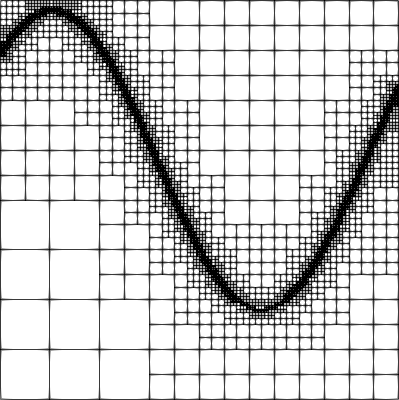}
    \end{tabular}
    \label{fig:sine-ref}
  }
  \hspace{2cm}
  \subfigure[``annulus'' refinement. Leaf mesh with finest cells on level 9.]
  {
    \begin{tabular}{ccc}
        \includegraphics[width=.3\textwidth]{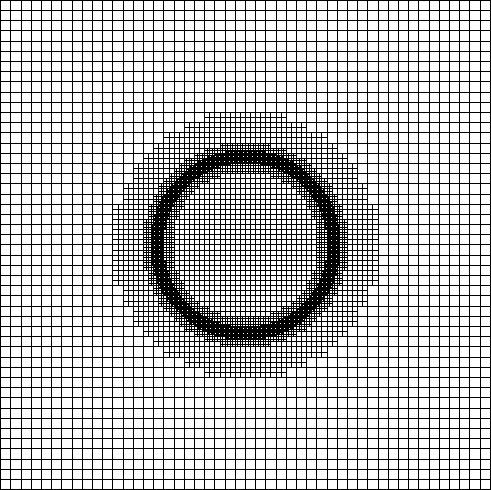}
    \end{tabular}
    \label{fig:sphere_shell-ref}
  }
  \caption{Visualization of the different mesh refinement sequences. See text for descriptions of each refinement algorithm.}
  \label{fig:sequences}
\end{figure}

Making general conclusions about the partitioning efficiency is
difficult as it depends on the number of processors, the coarse mesh, and
the refinement.  Instead, we study the efficiency for several test
cases shown in Figure~\ref{fig:sequences}.
These are obtained by the following construction. All are based on a coarse mesh consisting
of a single cell defined by $[-1,\,1]^{d}$. Finer meshes, where $L$ denotes the level of the finest cells, are obtained recursively by one of the following selection criteria:
\begin{description}
\item[``uniform''] global refinement of the coarse mesh $L$ times,
  obtaining a uniform leaf mesh of $4^L$ cells in two dimensions.
\item[``circle''] $L$ times refinement of all mesh cells with at least one
  vertex inside the circle of radius $\sfrac{1}{4\pi}$ around the origin.
\item[``quadrant''] $L$ times refinement of all mesh cells in the negative quadrant.
\item[``sine curve''] $L$ times refinement of all mesh cells intersected by the curve (plane for 3d) $\frac{3}{4}\sin\left(3x-\frac{5}{2}\right)+\frac{2}{10}$.
\item[``annulus''] After $L-3$ uniform refinements, add the steps:
  \begin{enumerate}
  \item Refine all cells whose center lies in the circle (sphere for 3D) of radius $0.55$.
  \item Refine all cells in the shell between radius $0.3$ and $0.43$.
  \item Refine all cells in the shell between radius $0.335$ and $0.39$.
  \end{enumerate}
\end{description}
All these procedures are completed by a closure after each refinement
step, ensuring one-irregularity in the sense that two leaf cells may only differ by one level, if they share a degree of freedom or a face for conforming methods and for DG methods, respectively.
These conditions are imposed in the {\dealii{}} library for practical
reasons, because they simplify several aspects of the
implementation. They can be relaxed at the price of software complexity.

\definecolor{color1}{HTML}{0A5F02}
\definecolor{color2}{HTML}{1C567A}
\definecolor{color3}{HTML}{814292}
\definecolor{color4}{HTML}{D7383B}
\definecolor{color5}{HTML}{9A7611}
\definecolor{color6}{HTML}{44AF7A}

\begin{figure}[tp]
  \begin{center}
    \footnotesize
    \tikzsetnextfilename{tikz_ratio/circle-quad}
    \includegraphics[width=0.49\textwidth]{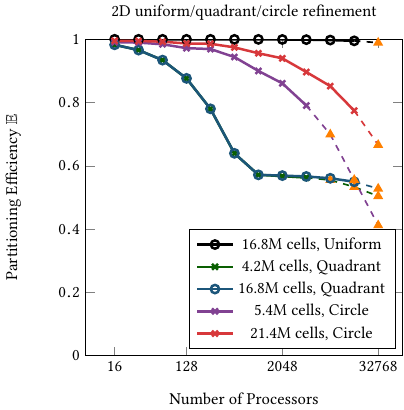}
    \hfill
    \tikzsetnextfilename{tikz_ratio/sine-curve}
    \includegraphics[width=0.49\textwidth]{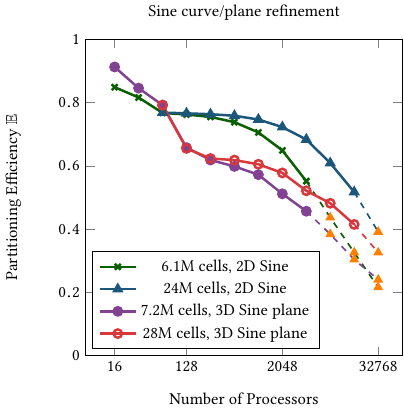}\\

    \tikzsetnextfilename{tikz_ratio/2D-martin_ref}
    \includegraphics[width=0.49\textwidth]{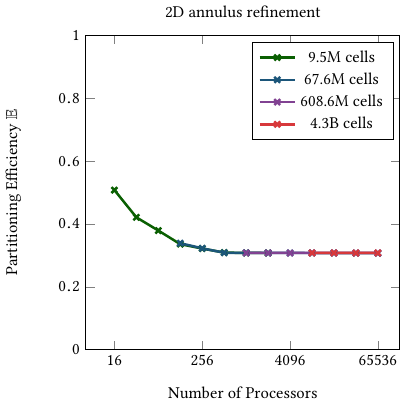}
    \hfill
    \tikzsetnextfilename{tikz_ratio/3D-martin_ref}
    \includegraphics[width=0.49\textwidth]{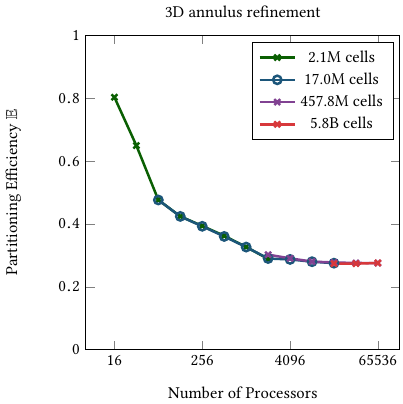}
  \end{center}
  \caption{Partitioning efficiency for various meshes. Dashed lines and orange triangles indicate experiments with less than 1,000 leaf cells per processor. Top: uniform, circle, quadrant refinement (left) and sine curve refinement (right). Bottom: annulus refinement in two (left) and three (right) dimensions.}
  \label{fig:ratio}
\end{figure}
Figure \ref{fig:ratio} shows the partitioning efficiency $\parteff$ as defined in \eqref{eq:slowdown-model} for varying processor count and problem size. For uniformly refined meshes we observe 100\% efficiency (this also holds in 3d, not shown). This is due to the fact that processor counts are multiples of $2^d$, which implies perfect partitioning on each level.
The ``quadrant'', ``sine curve'', and ``annulus'' refinement schemes show roughly the same behavior. Their efficiency drops until it levels off at 60\% for the quadrant and 2d sine curve, 50\% for the 3d sine plane, and about 30\% for the annulus in two and three dimensions. This level
of efficiency is then maintained over a wide range of processor counts. It only begins dropping again when the problem size is down to less than 1000 cells per processor, especially for the ``circle'' and ''sine curve'' refinements.
Both are extreme cases with very localized refinement and these examples have fewer cells compared to many of the other examples. For realistic problem sizes (more than 1000 active cells per processor), their efficiency remains above 50\%.

All of this information together suggests that, given a sufficient number of cells per processor, the imbalance of this distribution is primarily dependent on the type of mesh refinement refinement scheme. The number of processors influences the efficiency only to a certain ``leveling off'' point. In all cases, the efficiency stays above 30\% compared to the optimal
workload. In order to reach higher efficiency, re-partitioning and additional communication to processes farther away would be necessary.

Table~\ref{tab:workload_comm_levels} gives a level-by-level breakdown of the
partition efficiency for the 3D annulus refinement with 4.1M cells on 1,024
processors. Predictably, the efficiency is very low at lower levels where there are
fewer cells, and trends towards 1 for higher levels. It should be noted that the partition
efficiency calculated using $W_{\text{sync}}$ would be very different on lower levels than the
value $\parteff$ in the table which uses $W_{\text{opt}}$ (this is seen on the smaller
example in Section~\ref{sec:partition-efficiency}). However, the difference in the
total efficiency over all levels is only around 0.1\% with a value of 0.309 using $W_{\text{sync}}$ as
compared to 0.308 using $W_{\text{opt}}$.

\begin{table}[tp]
\caption{Workload imbalance and communication ratio per level for 3D annulus refinement,
4.1M cells, 1,024 processors. For partition efficiency, $W$, $W_{\text{sync}}$, $W_{\text{opt}}$,
and $\parteff$ are defined in Section~\ref{sec:partition-efficiency}. As discussed there, there
partition efficiency $\parteff$ is calculated using $W_{\text{opt}}$ instead of $W_{\text{sync}}$.
For communication ratio, ``ghost'' represents the number of cells which require communication
during the transfer from that level and ``native'' represents the number of cells which do not
require communication.}
\label{tab:workload_comm_levels}
 \begin{tabular}{|r|r||r|r|r|r||r|r|r|} \hline
   \multicolumn{2}{|c||}{}  & \multicolumn{4}{c||}{Partition efficiency} & \multicolumn{3}{c|}{Communication ratio}\\
   \hline
  level & cells & $W$ & $W_{\text{sync}}$ & $W_{\text{opt}}$ & $\parteff$ & ghost & native & ratio \\
  \hline
  0 & 125 & 1 & 1 & 0.1221 & 0.12207 & 513 & 487 & 0.51300 \\
  1 & 1,000 & 8 & 1 & 0.9766 & 0.12207 & 866 & 7,134 & 0.10825 \\
  2 & 8,000 & 64 & 8 & 7.8125 & 0.12207 & 2,028 & 61,972 & 0.03169 \\
  3 & 64,000 & 506 & 63 & 62.500 & 0.12352 & 3,118 & 508,882 & 0.00609 \\
  4 & 512,000 & 4,048 & 500 & 500.00 & 0.12352 & 3,017 & 354,743 & 0.00843 \\
  5 & 357,760 & 3,988 & 350 & 349.38 & 0.08761 & 2,635 & 807,349 & 0.00325 \\
  6 & 809,984 & 2,316 & 791 & 791.00 & 0.34154 & 0 & 2,977,280 & 0.00000  \\
  7 & 2,977,280 & 4,048 & 2,908 & 2,907.5 & 0.71826 & - & - & -  \\
  \hline
  Total & 4,730,149 & 14,979 & 4,622 & 4,619.3 & 0.30838 & 12,177 & 4,717,847 & 0.00257  \\
  \hline
\end{tabular}
\end{table}

\subsection{Communication}

The second factor determining the performance of parallel algorithms,
next to load balancing discussed above, is communication
overhead. Communication is not only much slower than computation at the granularity of individual instructions, it
also consumes more energy, because electrical charges must be
transported over fairly long distances. We introduced the first-child
rule with the express purpose to reduce communication overhead and keep
it small compared to the local computations. In
this section, we set out to demonstrate that this goal was achieved.

Communication happens in matrix-vector products and in grid transfer operations. Both of them apply a linear operator to a global discretization vector. The communication overhead in the first case is reduced by partitioning the leaf mesh into subdomains, such that their surface per volume ratio is small. Since the surface is of lower dimension, this implies that communication cost tends to zero as the number of cells on each processor grows to infinity. For weak scaling, this implies that it remains small compared to local operations, as long as there are sufficiently many cells on each processor. Such a partitioning is efficiently achieved by a space filling curve, in our case, the $z$-curve. Enumerating the cells along such a curve implies that cells with close indices will typically be close geometrically. This approach has been a standard for many years now.

\begin{figure}[tp]
	\begin{center}
		\footnotesize
		\tikzsetnextfilename{tikz_ratio/circle-quad-comm}
		\includegraphics[width=0.49\textwidth,height=0.35\textheight]{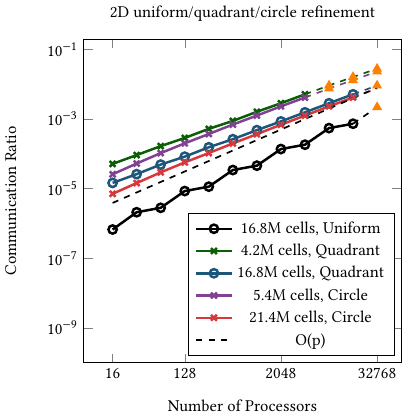}
		\hfill
		\tikzsetnextfilename{tikz_ratio/sine-curve-comm}
		\includegraphics[width=0.49\textwidth,height=0.35\textheight]{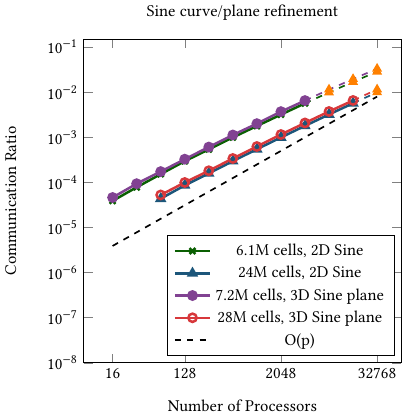}\\

		\tikzsetnextfilename{tikz_ratio/2D-martin_ref-comm}
		\includegraphics[width=0.49\textwidth,height=0.35\textheight]{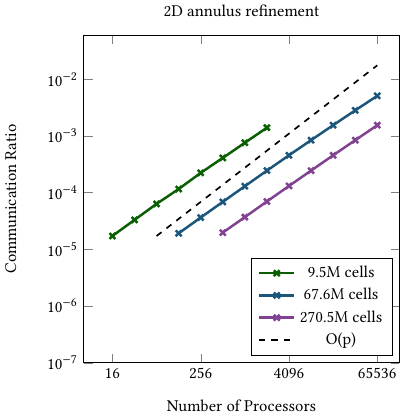}
		\hfill
		\tikzsetnextfilename{tikz_ratio/3D-martin_ref-comm}
		\includegraphics[width=0.49\textwidth,height=0.35\textheight]{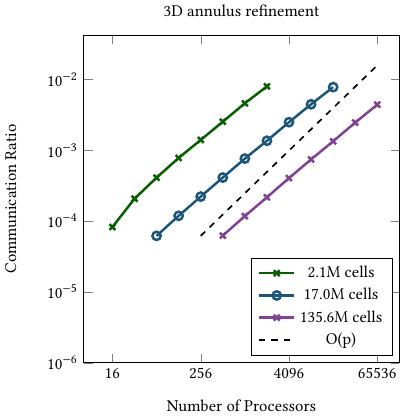}
	\end{center}
        \caption{Communication ratio (communicated number of children over total number) in level transfer
          for different test problems. The black
          dashed line denotes linear growth with number of processes, $O(p)$. Lines turn into dashed lines
          with orange triangles for data points with less than 1,000 leaf cells per processor. The number of
        cells to communicate grows more slowly than linear and stays below 1\%.}
  \label{fig:communication-ratio}
\end{figure}

The first child rule for partitioning lower levels
achieves a similar goal for grid transfer operations. If most
children are on the same processor as their parents, the amount of
communicated data is also much lower than the total amount of data
processed. In Figure~\ref{fig:communication-ratio}, we show that the
number of ``ghost children'' is indeed very small compared to the total
number of children. And while these numbers are rising with the number of processors, in the worst case observed less than 1\% of the cells require communication.
Additionally, the total communication volume seems to grow more slowly than the number of processors involved in the communication.

Finally, Table~\ref{tab:workload_comm_levels} gives a level-by-level
breakdown of the communication ratio for the 3D annulus refinement with
4.1M cells on 1,024 processors. Each of the higher levels require the communication
of about 2,000-3,000 cell during the transfer, with
the exception of the highest level, which requires no communication since
the p4est distribution used for the active level requires that all groups of
terminal children are owned by the same processor. For refinement schemes
where the majority of active cells are of the
highest refinement level, this will result in low communication ratios, as is seen with global refinement in Figure~\ref{fig:communication-ratio}.

\section{Performance Results}\label{sec:performance}

The algorithm described in the previous sections has been implemented in
the {\dealii{}} finite element library \cite{BangerthHartmannKanschat2007,dealII91}. The partitioning of
the adaptively refined meshes uses
p4est \cite{BursteddeWilcoxGhattas11}. The implementation with sparse
matrices uses Trilinos EPetra~\cite{trilinos-overview}, while the matrix-free
implementation is based on data distribution algorithms built into \dealii{}.
The source code and parameters of the examples in this manuscript are available at
\url{https://github.com/tjhei/paper-parallel-gmg-data}.

\subsection{Scaling on SuperMUC}

As a first experiment, we consider the constant-coefficient Laplacian on a cube,
discretized with $\mathcal Q_2$ elements, and compare the runtime on a uniform
mesh against an adaptively refined case with the annulus refinement. The
adaptive mesh is set up such that the number of cells matches with the number
of cells in the uniform case within 2\%. The computations are run on phase 1
of SuperMUC, providing nodes with $2\times 8$ cores of Intel Xeon E5-2680
(Sandy Bridge), connected via an Infiniband FDR10 fabric. For pre- and
post-smoothing, a Chebyshev iteration of the Jacobi method with Chebyshev
degree five, i.e., five matrix-vector products, is selected
\cite{adams03smoothing}. The relatively high degree of the Chebyshev smoother is
the result of an experimental study for the best run time for the chosen
implementation described below on uniform grids. Most important to the decision is the
mixed-precision setup with smoothing done in single precision as described in
\cite{kronbichler2019gpu,fehn2019hybrid}.
The parameters of the Chebyshev polynomial are set to
damp contributions in the eigenvalue range $[0.08 \bar{\lambda}_{\text{max},\ell}, 1.2 \bar{\lambda}_{\text{max},\ell}]$
on each level $\ell >0$. The estimate  $\bar{\lambda}_{\text{max},\ell}$ of the
largest eigenvalue of the matrix $A_\ell$ is computed by a conjugate gradient
iteration with 10 iterations from an initial vector of zero mean constructed
as $(-5.5, -4.5, \ldots, 4.5, 5.5, -5.5, -4.5, \ldots)^T$. As a coarse solver,
the Chebyshev iteration is selected with a degree chosen such that a priori
error estimate of the Chebyshev iteration ensures a residual reduction by
$10^3$, now for the full eigenvalue range of the coarse level matrix
determined by a conjugate gradient solution to a relative tolerance in the
unpreconditioned residual of $10^{-3}$. The coarse grid tolerance is chosen as
a balance between the accuracy needed to not influence the convergence rate
of the V-cycle for the chosen smoother (see section 5.1 of~\cite{kronbichler2019gpu}
for a detailed study) and its cost. Sine the coarse grid only consists of a
single cell in this experiment, the number of iterations is less than 5 in both 2D and 3D.

In order to reveal possible communication bottlenecks, we choose a fast node-level implementation by matrix-free evaluation of the matrix-vector products both for level matrices $A_\ell$ and level transfer \cite{kronbichler2016comparison}. The implementation exploits SIMD vectorization across several cells \cite{kronbichler2012generic} using four-wide registers on the given Intel Xeon processors. To further enhance performance, we run the multigrid V-cycle in single-precision as suggested originally in \cite{gropp00singleprec}. When combined with a correction in double precision after each V-cycle, e.g.~within an outer conjugate gradient solver, the reduced precision (which is of high-frequency character and thus easily damped in subsequent cycles) does not substantially alter the multigrid convergence \cite{kronbichler2019gpu}.

\begin{figure}[tp]
		\footnotesize
		\subfigure[2D. Left: 9.5M cells/37.9M DoFs. Right: 606.3M cells/2425.4M DoFs]{
		\tikzsetnextfilename{tikz_scaling/strong-2D-9M}
		\includegraphics[width=0.46\textwidth]{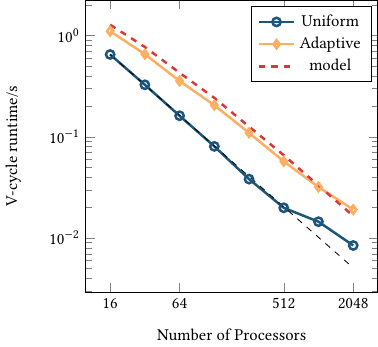}
		\hfill
		\tikzsetnextfilename{tikz_scaling/strong-2D-600M}
		\includegraphics[width=0.46\textwidth]{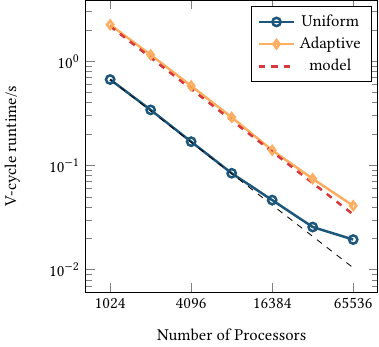}}

		\subfigure[3D. Left: 16.9M cells/137.4M DoFs. Right: 5.8B cells/46.4B DoFs]{
		\tikzsetnextfilename{tikz_scaling/strong-3D-17M}
		\includegraphics[width=0.46\textwidth]{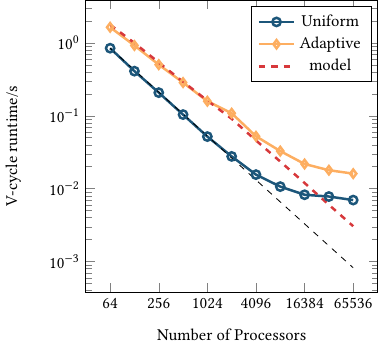}
		\hfill
		\tikzsetnextfilename{tikz_scaling/strong-3D-5B}
		\includegraphics[width=0.46\textwidth]{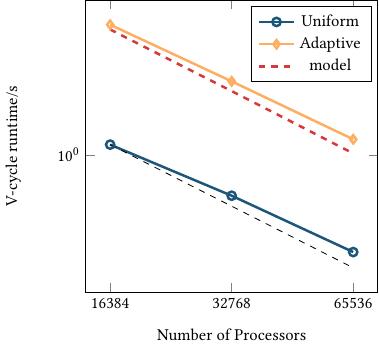}
		}
	\caption{Strong scaling for timing of a matrix-free V-cycle in 2D and 3D for small (left) and large (right) problem size of the annulus refinement.}
	\label{fig:strong}
\end{figure}

\begin{figure}[tp]
	\footnotesize
	\subfigure[2D annulus: 407K cells/1.6M DoFs per core]{
		\tikzsetnextfilename{tikz_scaling/weak-2D}
		\includegraphics[width=1.0\textwidth,height=0.32\textheight]{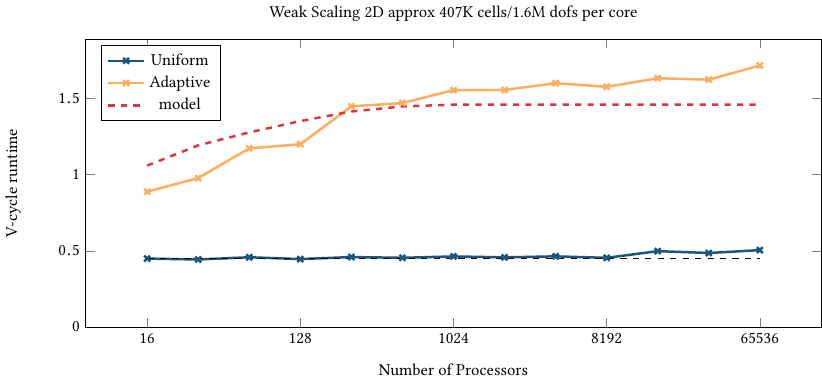}}

	\subfigure[3D annulus: 260K cells/2.1M DoFs per core]{
		\tikzsetnextfilename{tikz_scaling/weak-3D-reduced}
		\includegraphics[width=1.0\textwidth,height=0.32\textheight]{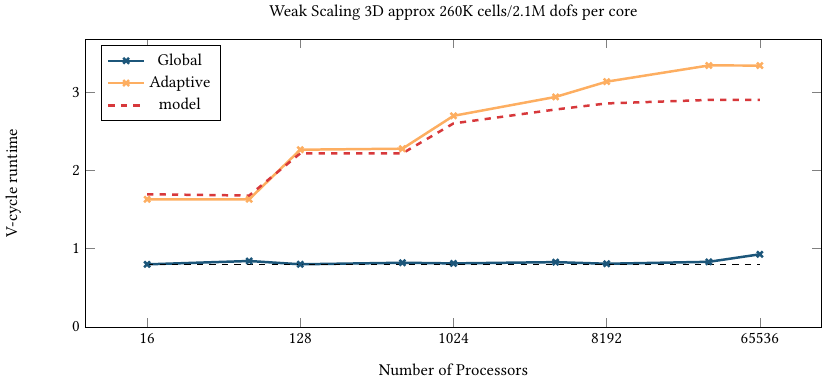}}
	\caption{Weak scaling for timing of a matrix-free V-cycle in 2D and 3D for the annulus refinement.}
	\label{fig:weak}
\end{figure}

Figures \ref{fig:strong} and \ref{fig:weak} list the strong and weak scaling for
the runtime of one multigrid V-cycle run as a preconditioner, including all
aforementioned communication steps as well as the conversion from double to
single precision and vice versa. The presented numbers are consistent over
several runs (with standard deviations of at most 2\% of the runtime). Each
plot contains runtimes for the uniform and adaptive refinement and the optimal
$\mathcal{O}(N)$ scaling (black dashed line) coinciding with the first data
point of the uniform refinement graph.
The red dashed line shows the model prediction based on the imbalance
of the adaptively refined mesh computed as $1/\parteff$ multiplied by the
ideal scaling of the uniform computation for the same processor counts (black dashed line).
Given the results in Figure~\ref{fig:ratio},
the 2D annulus refinement suggests an efficiency gap of a factor close to 3
in two dimensions.
Figure~\ref{fig:weak}(a) confirms this behavior, confirming that the model
assumption is realistic: the uniform refinement is predicted to be 100\%
efficient and the adaptive refinement is 31\% efficient, so we predict a
gap of 3.2$\times$ in runtime.

Finally, it is worth mentioning that we only looked at the performance of a single
V-cycle so far. When using it as a preconditioner of a Conjugate Gradient solver, the performance
results are very
similar as the number of iteration stays nearly constant.
In fact, the number of iterations is between 7 and 11 for a relative reduction of $1e-10$
for uniform and adaptive runs in 2d and 3d.


\definecolor{color1}{HTML}{1C567A}
\definecolor{color2}{HTML}{FDAE61}
\definecolor{color3}{HTML}{D7383B}

The strong scaling limit of the adaptive implementation follows the one of the
uniform case, highlighting the efficiency in the communication setup. In three
dimensions with 16.9 million cells, scaling of the uniform mesh case starts to
flatten for 8,192 MPI ranks, corresponding to 2048 cells or approximately
54,000 unknowns (DoFs) per MPI rank. For this data point, the absolute runtime for the
V-cycle is 0.01 seconds. Given the fact that 12 matrix-vector products are
performed per level (8 in the smoother, two for the residuals, two for the
transfer) for a total of 8 levels, this data point corresponds to
approximately $1.0\cdot 10^{-4}$ seconds per matrix-vector product, which is
an expected scaling limit of nearest neighbor communication for up to 26
neighbors combined with some local computation on the given architecture.
The adaptive case scales at least as well as the uniform one even beyond 8k
processors, and also for the other experiments. Partly, this is due to an overlap
of different levels e.g. when some processors do not own any part of a fine
level, they can start working on coarser levels as long as the local
communication data arrives. Furthermore, the imbalance also leads to more
cells on the processors for a given level in relative terms approximately
proportional to the inverse efficiency factor $1/\parteff$.

\subsection{Linear elasticity with discontinuous Galerkin discretization}
\label{sec:test-elasticity}

\begin{figure}[tp]
  \begin{center}
    \footnotesize
    \includegraphics[width=0.48\textwidth]{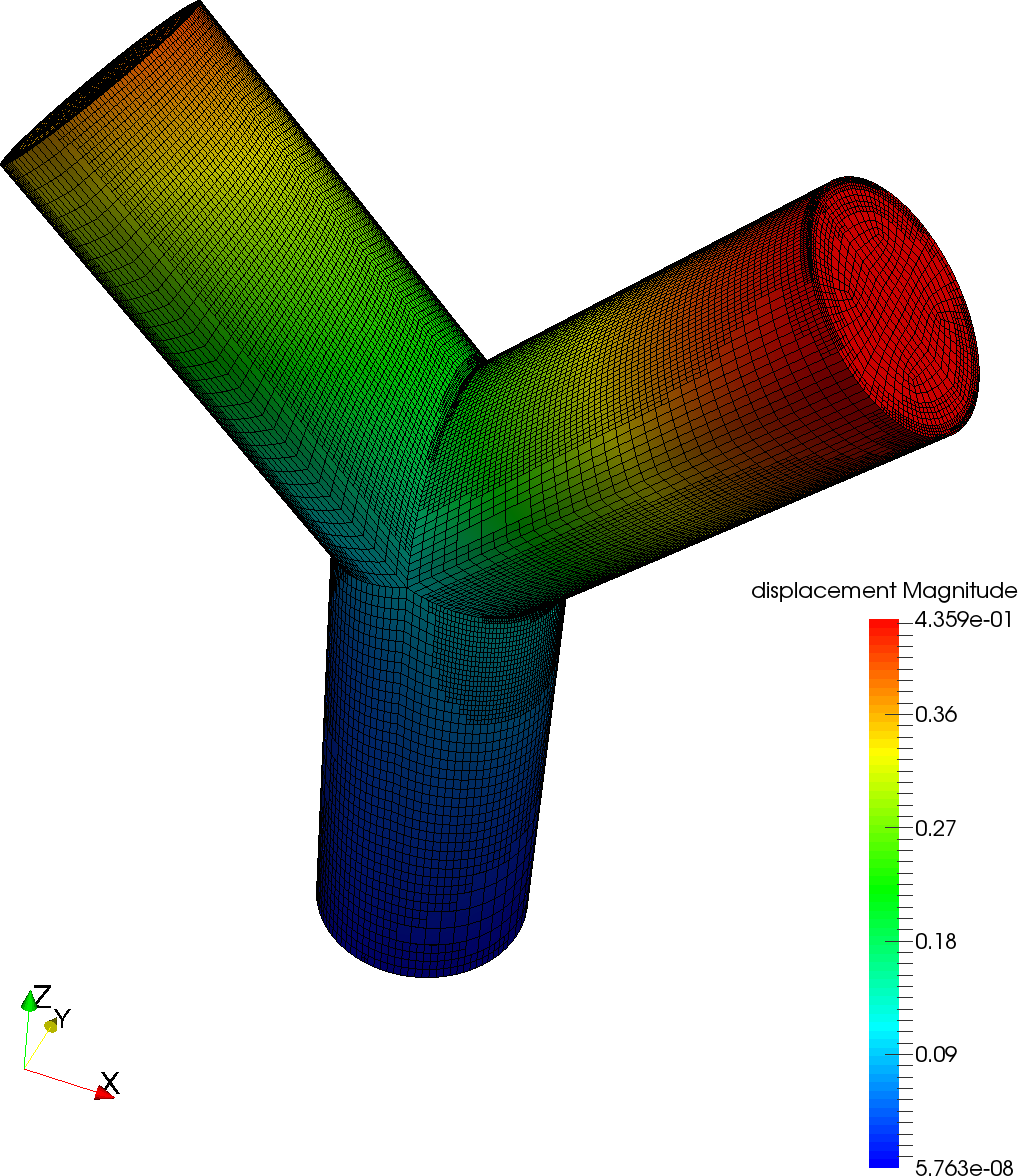}
    \hfill
    \includegraphics[width=0.48\textwidth]{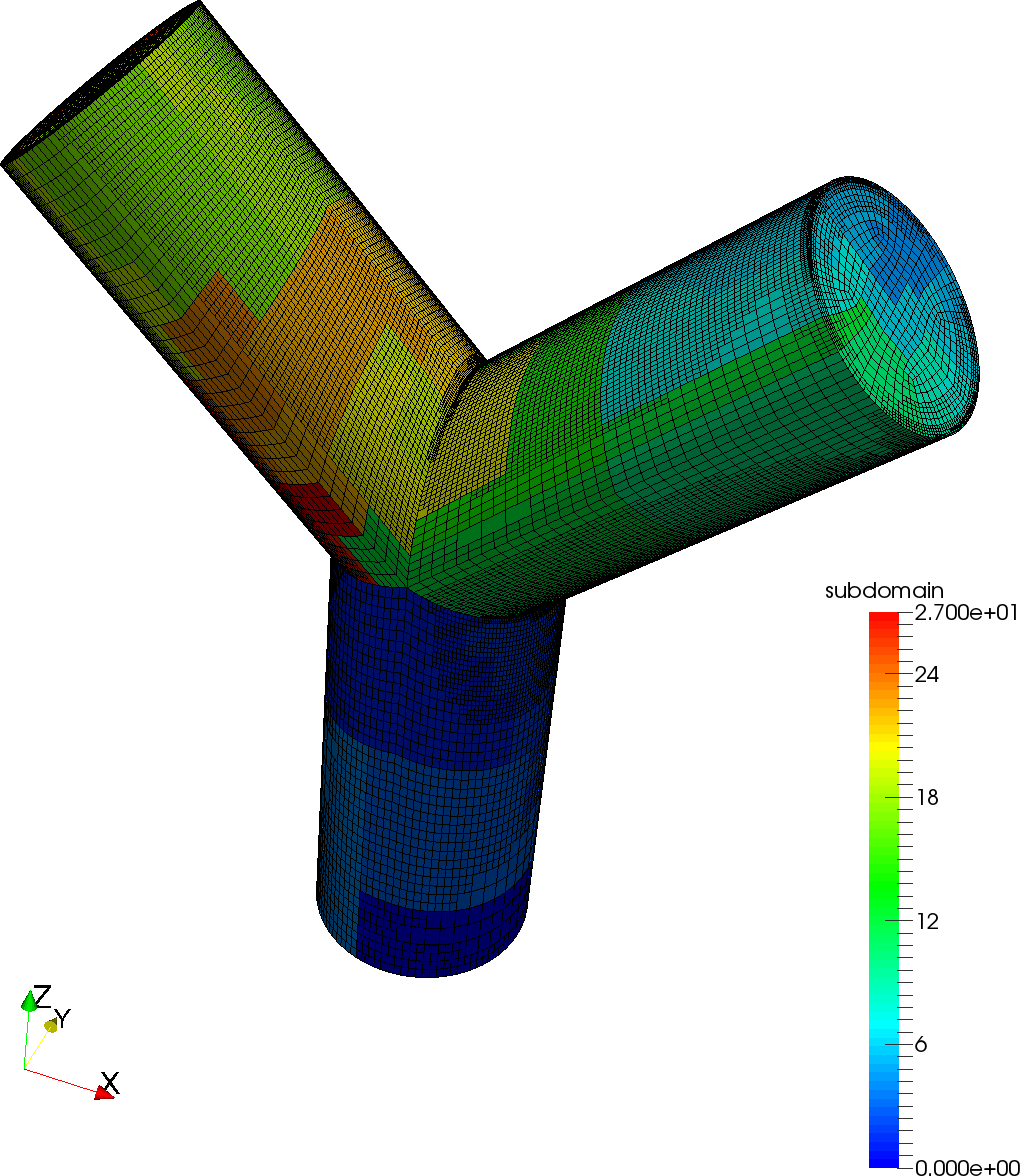}
  \end{center}
  \caption{Displacement magnitude (left) and distribution of domain over
    28 processors (right) of three-dimensional elasticity example.}
  \label{fig:elasticity}
\end{figure}

As a second experiment, we consider the equations of linear elasticity
\eqref{eq:elasticity} on a mesh constructed from three cylinders with the
Lam\'e parameters $\lambda = \mu = 1$ according to the setup in
Fig.~\ref{fig:elasticity}. The solid is loaded by surface forces on the upper bases of the top
two cylinders. It is fixed at the base of the lower cylinder and traction-free on the sides the cylinders.
In order to represent the geometry with a high-quality
mesh, we use 2808 hexahedral cells with one global and a series of up to
three adaptive refinements based on a residual-based error
estimator. Fig.~\ref{fig:elasticity} shows how the error estimator chooses to
refine around the sharp corners with lower solution regularity.
The outer layer of cells is represented by a curved cylindrical
manifold aligned with the respective cylinder sides. To smoothly relax the curved
surface description into a straight-sided one towards the center of
the cylinder, we apply a transfinite interpolation \cite{gordonthiel82} over
approximately half the cylinder radius. For approximation, we use
vector-valued discontinuous $\mathcal Q_2$ elements of tensor degree 2 and the
symmetric interior penalty method with penalty factor $\sigma_h$
equal to 2.0 weighted by the
minimum vertex difference in face-normal direction and the factor
$2\cdot 3 =6$ to account for the inverse estimate on quadratic shape
functions.

We solve the elasticity example with a point-Jacobi smoother with four pre-
and postsmoothing sweeps and relaxation parameter 0.5 on all levels, using a
matrix-based implementation based on Trilinos EPetra linear algebra.
On the coarse level, there are 227K ($=2808\times 81$)
unknowns and 123M nonzero entries in the matrix. We compare two different
strategies for solving this coarse linear system. The first setup uses a
direct solver based on the
\mbox{SuperLUDist} package, whereas the second uses an iterative conjugate gradient
solver preconditioned by the Trilinos AMG preconditioner ML. The coarse grid CG
solver is run to a relative tolerance of $10^{-2}$, compared against the initial
unpreconditioned residual. The AMG solver is given the near-null space of
elasticity, i.e., three translational and three rotational modes, the latter
using the coordinates of the nodal points of the finite element interpolation. Two
sweeps of an incomplete LU factorization  (no fill-in, no overlap in parallel) are used
for pre- and post-smoothing, and the aggregation threshold is 0.01.
All other setting use the standard settings for
elliptic problems.

The systems are then solved by a conjugate gradient solver on the leaf mesh
preconditioned by the proposed geometric multigrid scheme to a relative
tolerance of $10^{-6}$, measured in the unpreconditioned residual norm.
Table~\ref{tab:elasticity} displays the number of
iterations and runtimes on 28 processors (one node) for the two options. The results demonstrate that
the multigrid preconditioner yields mesh-independent iteration counts also for the
elasticity problem and a more complex geometry. In particular, the run time per unknown
is constant or even slightly decreases as the grid is refined, showing that all
components in the multigrid algorithm show optimal weak scaling as the problem
size is increased. However, the iterative
coarse-grid solver produces solver runtimes which are considerably worse than
the direct solver SuperLUDist. The high cost of the iterative solver is due
to the large number of iterations. For the example of 11.7 million unknowns,
the coarse solver takes 27.8 iterations for each outer CG iteration on
average (or 779 when accumulating over all iterations). This high
iteration count is due to the higher-order discontinuous nature of the
solution space and could be overcome, e.g., by p-multigrid techniques \cite{helenbrook03},
see also \cite{fehn2019hybrid} and references therein.

\begin{table}
\caption{Number of outer conjugate gradient iterations and runtimes for solving the elasticity example on a 28-core setup for two strategies on the coarse mesh.}
\label{tab:elasticity}
 \begin{tabular}{|r|r|r|r|r|r|r|r|r|} \hline
    &  & \multicolumn{3}{c|}{coarse direct solver} & \multicolumn{4}{c|}{coarse CG/AMG iteration}\\
   & & &  time &  time / DoF & &  coarse solve &  time  & time / DoF\\
   levels & DoFs & CG its & [s] & [$\mu$s / DoF] & CG its & avg CG its & [s] &  [$\mu$s / DoF]\\ \hline
   2 & 1,819,584 & 27 & 46.1 & 25.3 & 27 & 28.8 & 113.4 & 62.3 \\
   3 & 2,456,325 & 27 & 63.6 & 25.9 & 27 & 29.7 & 152.7 & 62.2 \\
   4 & 4,916,538 & 28 & 121.7 & 24.8 & 28 & 29.7 & 251.9 & 51.2 \\
   5 & 11,684,817 & 28 & 260.3 & 22.3 & 28 & 27.8 & 420.0 & 35.9 \\ \hline
  \end{tabular}
\end{table}

Table~\ref{tab:elasticity_scaling} shows a scaling experiment of the
elasticity example on up to 3,072 processors with the coarse direct solver
SuperLUDist. In order to separate the lack of parallel scaling in the direct
solver from the proposed contributions, the table also reports the solver time
for all parts except the coarse grid solver, derived from subtracting timing
results of the coarse grid solver with an MPI barrier around it from the
overall solver time. For all cases, the number of
solver iterations is between 27 and 29, similar to the results reported in
Table~\ref{tab:elasticity}. The results demonstrate the ability of the
proposed framework to run on a large scale also when using matrix-based
solvers. However, due to the discontinuous $\mathcal Q_2$ basis, the average
number of nonzero entries per row exceeds 500, leading to a high memory
consumption with up to 4 GB per MPI process of resident virtual memory for the
largest computation of 317 million unknowns on 3,072 cores. Larger problems
for a given memory configuration could be solved with matrix-free
solvers~\cite{kronbichler2019fast}.

\begin{table}
\caption{Runtimes in seconds for solving the elasticity example on a system of 48-core nodes for a coarse direct solver.}
\label{tab:elasticity_scaling}
 \begin{tabular}{|r|r|r|r|r|r|r|r|r|r|} \hline
   &  & \multicolumn{2}{c|}{1 node} & \multicolumn{2}{c|}{4 nodes} &  \multicolumn{2}{c|}{16 nodes} &  \multicolumn{2}{c|}{64 nodes} \\
   & & & without & & without & & without & & without \\
   levels & DoFs & time & coarse & time & coarse & time & coarse & time & coarse\\
   \hline
   2 & 1.8M    & 27.3  & 14.9  & 35.5  & 5.2   & 30.3  & 1.5   & 40.9  & 0.4   \\
   3 & 2.5M    & 44.5  & 28.7  & 32.7  & 8.4   & 30.5  & 3.0   & 30.1  & 1.1   \\
   4 & 4.9M    & 85.9  & 71.5  & 51.2  & 21.4  & 34.9  & 7.8   & 45.8  & 2.8   \\
   5 & 11.7M   & 188.3 & 173.1 & 82.7  & 54.4  & 48.9  & 18.2  & 35.1  & 6.3   \\
   6 & 35.7M   & ---   & ---   & 212.1 & 185.3 & 91.3  & 62.1  & 48.8  & 20.6  \\
   7 & 101.1M  & ---   & ---   & ---   & ---   & 214.2 & 178.9 & 92.8  & 61.6  \\
   8 & 316.8M  & ---   & ---   & ---   & ---   & ---   & ---   & 262.3 & 227.8 \\
   \hline
  \end{tabular}
\end{table}


\begin{figure}[tp]
  \footnotesize
  \tikzsetnextfilename{tikz_scaling/elasticity-strong}
  \includegraphics[height=0.45\textwidth]{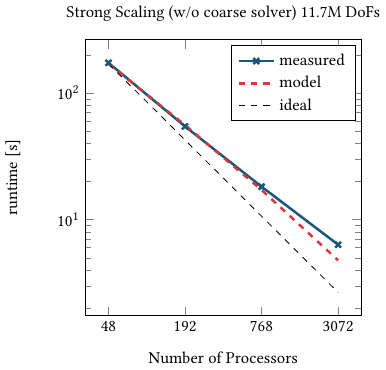}
  \hfill
  \tikzsetnextfilename{tikz_scaling/elasticity-weak}
  \includegraphics[height=0.45\textwidth]{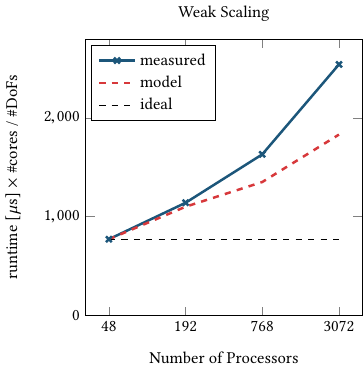}
  \caption{Parallel scaling for the elasticity problem in 3D. Strong scaling
    (left) uses a problem size of 11.7M DoFs and reports timings without the
    coarse solver. Weak scaling (right) is reported for the full solver times
    of the largest problem sizes per processor configuration reported in
    Table~\ref{tab:elasticity_scaling}, re-scaled to a fixed number of
    unknowns per MPI process.}
  \label{fig:elasticity_scaling}
\end{figure}

Figure~\ref{fig:elasticity_scaling} shows the strong and weak scaling behavior
of the data presented in Table~\ref{tab:elasticity_scaling} alongside the
efficiency model from Section~\ref{sec:partition-efficiency}. In the strong
scaling experiment on 11.7 million DoFs, the partitioning efficiency
$\mathbb{E}$ increases from 2.12 on 48 MPI ranks to 3.81 on 3,072 MPI
ranks. For weak scaling, the largest problem with 317 million DoFs on 8 levels
has an imbalance of 4.85. The model explains the loss of efficiency especially
for the strong scaling test and moderate sizes. For weak scaling, the measured
run time scales somewhat worse than the one predicted by the model, which can
be traced back to a memory bandwidth effect: On low node counts, the imbalance
is within a shared memory region, and cores with more work on a particular
level can use the memory bandwidth of idle cores. This effect mostly disappears for
higher node counts.

\subsection{Laplace: comparison against AMG}
\label{sec:test-laplace}

\begin{figure}[tp]
  \begin{center}
    \footnotesize
    \includegraphics[width=0.48\textwidth]{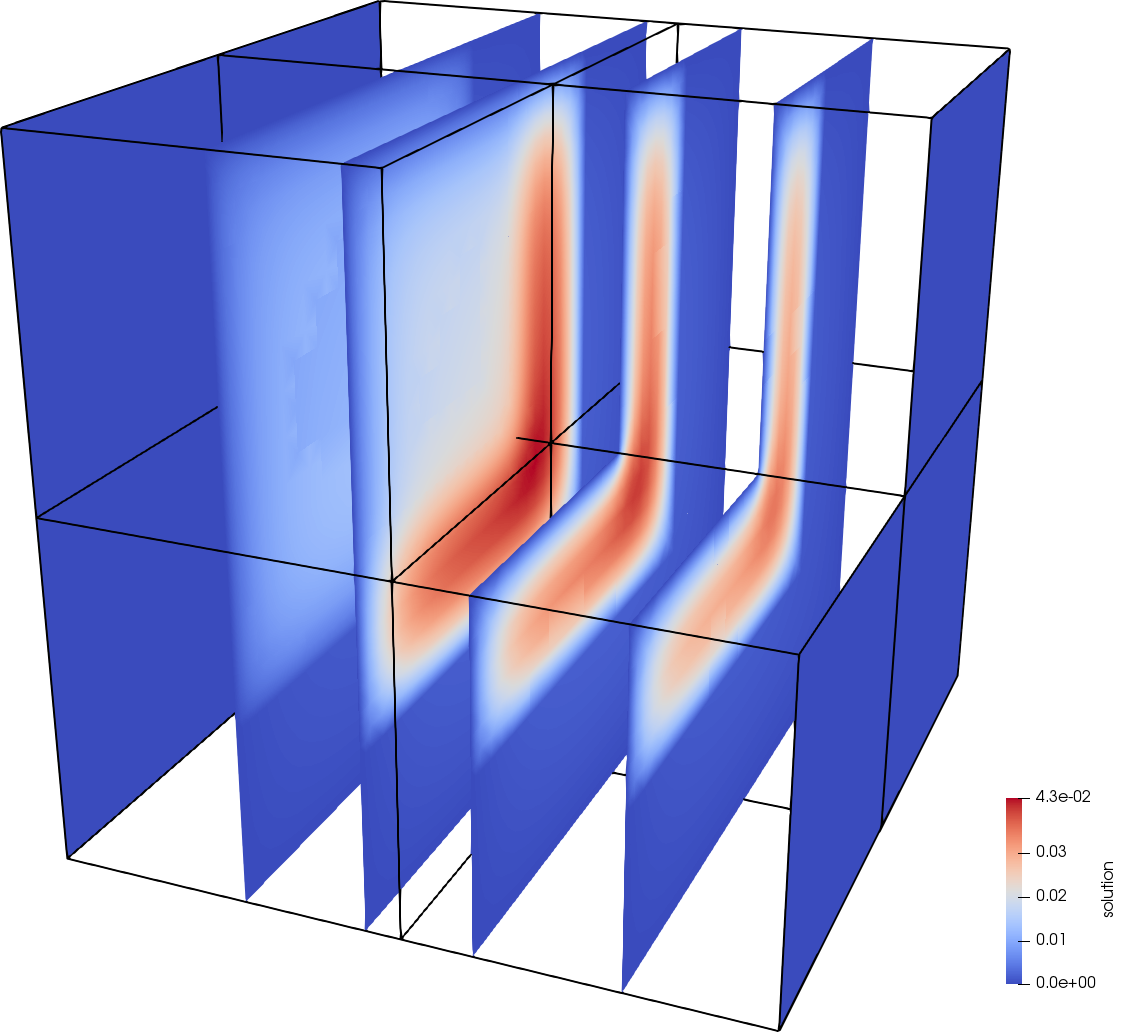}
    \hfill
    \includegraphics[width=0.37\textwidth]{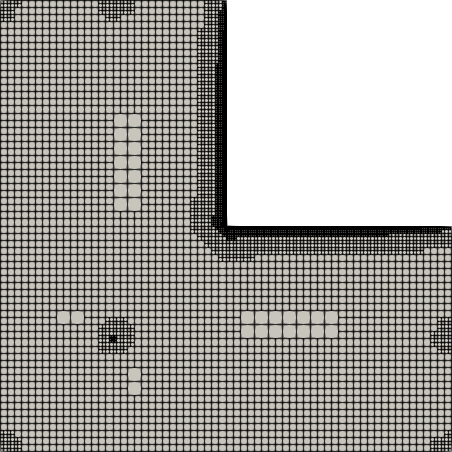}
  \end{center}
  \caption{Left: Solution to the Laplace example on a 3D domain with Fichera corner and discontinuous coefficients.
  Right: Slice for x close to the center of the domain showing the adaptively refined mesh.}
  \label{fig:laplace-3d}
\end{figure}

In this example we consider the variable-coefficient Laplacian
\[
 (\epsilon \nabla u, \nabla v) = (f,v) \quad \forall v \in V_h
\]
on the domain $\Omega = [-1,1]^3 \setminus [0,1]^3$ (a 3D Fichera corner) with
$\epsilon = 1$ if $\min(x,y,z)>-\frac{1}{2}$ and $ \epsilon = 100$ otherwise.
The boundary conditions are $u=0$ on the whole boundary and the right-hand side is $f=1$.
Figure~\ref{fig:laplace-3d} visualizes the solution.
We use continuous Q2 elements to discretize $V_h$ and use an adaptively refined grid.

For this example, we
compare a matrix-free geometric multigrid
implementation, a matrix-based geometric multigrid implementation,
and an algebraic multigrid based on Trilinos ML.
With this problem and this
discretization, we expect the algebraic multigrid method to work very well.
For a balanced comparison, we are picking the same settings
for all solvers, namely a
single Jacobi smoothing step. These might not be the optimal settings
for each individual
method, though. We compare the solution with the conjugate gradient method to a relative tolerance of $10^{-10}$ in the unpreconditioned residual norm with the
following preconditioners:
\begin{enumerate}
 \item ``MF'': matrix-free geometric multigrid using
 1 Jacobi step (without damping) in single precision, constructed all the way down to the coarsest
 mesh. The coarse solver is an unpreconditioned CG solve.
 \item ``MB'': geometric multigrid using Trilinos Epetra matrices using 1 Jacobi step (without damping), constructed all the way down to the coarsest
 mesh. The coarse solver
 is an unpreconditioned CG solve.
 \item ``AMG'': Trilinos ML using Epetra with
    1 Jacobi smoother, aggregation threshold 0.02, and coarse solve ``Amesos-KLU''.
\end{enumerate}

We start with a coarse mesh with 7 cells (and 117 Q2 unknowns) and refine
adaptively using the residual-based, cell-wise a posteriori error estimator $e(K) = e_{\text{cell}}(K) + e_{\text{face}}(K)$ from~\cite{karakashian2003posteriori} with
\[
 e_{\text{cell}}(K) = h^2 \| f + \epsilon \triangle u \|_K^2, \qquad
 e_{\text{face}}(K) = \sum_F h_F \| [ \epsilon \nabla u \cdot n ] \|_F^2.
\]
We double the number of cells approximately in each step of the
adaptive loop by refining the 14.2\% of the cells with the largest contribution to the estimator.

In Table~\ref{tab:laplace-example-time} we present the time in seconds
for the setup and for solving the linear system for the three methods and different problem sizes.
The different columns contain the following operations:
\begin{itemize}
\item
``mesh'': refine the mesh adaptively, repartition the cells between processors,
computing cell ownership, exchange ghost information (MF and MB).
\item
``FE'': enumerate the degrees of freedom, create vectors, create system matrix (MB, AMG).
\item
``assembly'': assemble system matrix (MB, AMG) and right-hand side, evaluate
viscosity (MF).
\item
``prec'': enumerate DoFs on each multigrid level (MB, MF), assemble level matrices (MB), create intergrid transfer operators (MB and MF), create matrix-free operators (MF), create AMG preconditioner (AMG).
\end{itemize}

The efficiency $\parteff$ for the four different problem sizes is
0.371, 0.294, 0.229, and 0.161, respectively. This can explain a slow-down
of $0.371/0.161 \approx 2.3$ the preconditioner setup from the smallest
to largest problem, which is close to what we observe here.

\begin{figure}[tp]
		\footnotesize
		\tikzsetnextfilename{tikz_laplace/strong}
		\includegraphics[width=0.76\textwidth,height=0.36\textheight]{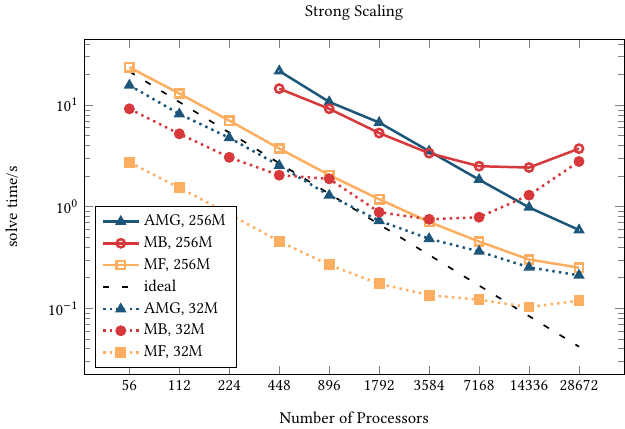}
\caption{Strong scaling for the solve step of the variable viscosity Laplace problem with adaptive refinement for two different meshes for algebraic multigrid (AMG), matrix based geometric multigrid (MB), and matrix-free geometric multigrid (MF).}
  \label{fig:laplace-strong}
\end{figure}

Finally, we present strong scaling of the CG solve for this highly
resolved mesh in Figure~\ref{fig:laplace-strong}.
While the matrix-based solver and algebraic multigrid scale similarly and
have a similar time to solution, the matrix-free implementation scales
much better and solves the finer problem in roughly the same time as
the AMG method for the coarser mesh with only an eighth of the number of unknowns.
For large numbers of processors, the solve time for the matrix-based solver
degrades considerably. Based on microbenchmarks, we hypothesize
that this effect is mostly caused by the implementation of the ghost exchange
of the \texttt{Epetra\_CrsMatrix} in Trilinos, which includes a global barrier
and MPI ready-sends.\footnote{The implementation is given at
  \url{https://github.com/trilinos/Trilinos/blob/67064dc0bd94754b56ded6ee2a0a09a2dcc45433/packages/epetra/src/Epetra_MpiDistributor.cpp\#L781-L938}}
This operation involves all processors, including those idle on a given
level. For the high number of levels (16 for the 32M DoFs case, 19 for the
256M DoFs case) and the appearance both in the residual, the level transfer
computations, and the coarse grid CG solver, this operation affects the run
time. While the situation could be improved by suitable sub-communicators, the
much better scaling with the matrix-free solvers and their non-blocking MPI
communication show that this limit is not inherent to the proposed multigrid
algorithms but rather the linear algebra back-end.

\begin{table}[tp]
\caption{Time in seconds for various stages for different problem sizes for the three
solver schemes for the adaptively refined Laplace problem. MF: matrix-free geometric multigrid,
MB: matrix-based geometric multigrid, AMG: algebraic multigrid (Trilinos ML).}
\label{tab:laplace-example-time}
\begin{tabular}{|r|r|r||r|r|r|r|r|r|r|r|}
\hline
 &  &  &  & \multicolumn{5}{c|}{Setup}
 & & total\\
 & Proc & Cycle & DoFs & mesh & FE & assembly & prec & total & solve & time \\ \hline \hline
  MF & 112 & 13 & 4.1M & 0.562 & 0.151 & 0.029 & 0.393 & 1.135 & 0.200 & 1.335 \\ \hline
   & 448 & 15 & 16.3M & 0.703 & 0.154 & 0.027 & 0.535 & 1.419 & 0.253 & 1.672 \\ \hline
   & 1792 & 17 & 65.1M & 0.910 & 0.182 & 0.030 & 0.686 & 1.808 & 0.309 & 2.117 \\ \hline
   & 7168 & 19 & 256.3M & 1.030 & 0.152 & 0.032 & 0.893 & 2.107 & 0.521 & 2.628 \\ \hline \hline
  MB & 112 & 13 & 4.1M & 0.564 & 0.527 & 0.623 & 2.934 & 4.648 & 0.716 & 5.364 \\ \hline
   & 448 & 15 & 16.3M & 0.702 & 0.548 & 0.677 & 3.776 & 5.703 & 1.190 & 6.893 \\ \hline
   & 1792 & 17 & 65.1M & 0.898 & 0.575 & 0.698 & 4.862 & 7.033 & 1.660 & 8.693 \\ \hline
   & 7168 & 19 & 256.3M & 1.040 & 0.619 & 0.727 & 7.260 & 9.646 & 2.560 & 12.206 \\ \hline \hline
  AMG & 112 & 13 & 4.1M & 0.296 & 0.507 & 0.657 & 0.405 & 1.865 & 0.924 & 2.789 \\ \hline
   & 448 & 15 & 16.3M & 0.339 & 0.534 & 0.671 & 0.456 & 2.000 & 1.150 & 3.150 \\ \hline
   & 1792 & 17 & 65.1M & 0.421 & 0.553 & 0.680 & 0.546 & 2.200 & 1.460 & 3.660 \\ \hline
   & 7168 & 19 & 256.3M & 0.515 & 0.585 & 0.744 & 1.010 & 2.854 & 1.890 & 4.744 \\ \hline
  \end{tabular}
\end{table}

\section{Conclusions}\label{sec:conclusions}

In this article, we described the implementation of a parallel, adaptive multigrid framework within the multi-purpose finite element library \dealii{}. The framework allows for conforming as well as discontinuous finite elements on locally refined meshes. We have shown scaling results involving up to 65,536 cores with very good weak scaling and strong scaling as long as the local problem size is large enough.
The distribution of mesh hierarchies is optimized for communication reduction, such that the framework is expected to scale well after node-level optimizations through vectorization and algorithms with higher computational intensity. We exemplified the efficiency by evaluating the parallel scaling using a matrix-free implementation with optimized node-level performance.
We presented a model for the efficiency of the partitioning of the hierarchy and compared its prediction to actual runtimes. Computational experiments include an elastic structure with a nontrivial coarse mesh and comparison to algebraic multigrid.
The presented ingredients are flexible in terms of finite element spaces, matrix-based or matrix-free implementations, and smoothers.

The proposed performance model and the computational results suggest that
minimizing solely for communication is not optimal for performance. It is
subject of future research to identify the best tradeoff between load
imbalance and longer-distance communication in the level transfer in terms of
performance.

\section*{Acknowledgments}
Thomas C. Clevenger and Timo Heister
were partially supported by NSF Award OAC-2015848 and by the Computational Infrastructure in
Geodynamics initiative (CIG), through the NSF under Award EAR-0949446 and
EAR-1550901 and The University of California -- Davis.
Timo Heister was also partially supported by NSF Award
DMS-2028346, EAR-1925575, and by Technical Data
Analysis, Inc. through US Navy SBIR N68335-18-C-0011.
The work of Guido Kanschat and
Martin Kronbichler was supported by the German Research Foundation (DFG) via
the project ``High-order discontinuous Galerkin for the exa-scale''
(ExaDG) within the priority program 1648 ``Software for Exascale Computing''
(SPPEXA), grant agreements KA\,1304/2-1 and KR\,4661/2-1. Furthermore,
the support by the state of Baden-W\"urttemberg through bwHPC and the DFG
through grant INST 35/1134-1 FUGG are acknowledged.

The Gauss Centre for Supercomputing e.V. (www.gauss-centre.eu) funded this
project by providing computing time on the GCS Supercomputer SuperMUC at
Leibniz Supercomputing Centre (www.lrz.de) through project pr83te.  Clemson
University is acknowledged for generous allotment of compute time on Palmetto
cluster.
This work used the Extreme Science and Engineering Discovery Environment (XSEDE), which is supported by National Science Foundation grant number ACI-1548562.
The authors acknowledge the Texas
Advanced Computing Center (TACC) at The University of Texas at Austin for providing access to both the
Stampede2 and Frontera machines that have contributed to the research results reported within this paper.

The authors would like to thank their coauthors on the {\dealii{}} project.

\bibliographystyle{plain}
\bibliography{compiled,bib}
\end{document}